\documentclass[11pt,reqno]{amsart}
\usepackage{amsmath}
\usepackage{amssymb}
\usepackage{amsthm}

\newtheorem{theorem}{Theorem}[section]
\newtheorem{prop}[theorem]{Proposition}
\newtheorem{lem}[theorem]{Lemma}
\newtheorem*{cor}{Corollary}
\theoremstyle{definition}
\newtheorem{defn}[theorem]{Definition}
\theoremstyle{remark}
\newtheorem*{rem}{Remark}
\numberwithin{equation}{section}

\newcommand{\Gg}{\mathfrak{g}}    
\newcommand{\Gh}{\mathfrak{h}}
\newcommand{\Gq}{\mathfrak{q}}

\newcommand{\Gz}{\mathfrak{z}}

\begin{document}

\title{Haar measure on a locally compact quantum group}

\author{Byung-Jay Kahng}
\date{}
\address{Department of Mathematics\\ University of Kansas\\
Lawrence, KS 66045}
\email{bjkahng@math.ku.edu}
\thanks{The author wishes to thank Professor George Elliott, for his
kind and helpful comments on the early draft of this paper.}
\subjclass[2000]{46L65, 46L51, 81R50}
\keywords{locally compact quantum group, Haar weight, multiplicative unitary
operator}

\begin{abstract}
In the general theory of locally compact quantum groups, the notion of
Haar measure (Haar weight) plays the most significant role.  The aim of
this paper is to carry out a careful analysis regarding Haar weight, in
relation to general theory, for the specific non-compact quantum group
$(A,\Delta)$ constructed earlier by the author.  In this way, one can
show that $(A,\Delta)$ is indeed a ``($C^*$-algebraic) locally compact
quantum group'' in the sense of the recently developed definition given
by Kustermans and Vaes.  Attention will be given to pointing out the
relationship between the original construction (obtained by deformation
quantization) and the structure maps suggested by general theory.
\end{abstract}
\maketitle

{\sc Introduction.}
According to the widely accepted paradigm (which goes back to Gelfand 
and Naimark in the 1940's and has been reaffirmed by Connes and his
non-commutative geometry program \cite{Cn}) that $C^*$-algebras are
quantized/non-commutative locally compact spaces, the $C^*$-algebra
framework is the most natural one in which to formulate a theory of
{\em locally compact quantum groups\/}.  There have been several examples
of $C^*$-algebraic quantum groups constructed, beginning with Woronowicz's
(compact) quantum $SU(2)$ group \cite{Wr1}.  The examples of non-compact
$C^*$-algebraic quantum groups have been rather scarce, but significant
progress has been made over the past decade.

 Among the examples of non-compact type is the Hopf $C^*$-algebra $(A,\Delta)$
constructed by the author \cite{BJKp2}.  The construction is done by the method
of deformation quantization, and the approach is a slight generalization of the
one used in Rieffel's example of a solvable quantum group \cite{Rf5}.  In fact,
$(A,\Delta)$ may be regarded as a ``quantized $C^*(H)$'' or a ``quantized $C_0(G)$'',
where $H$ is a Heisenberg-type Lie group and $G$ is a certain solvable Lie group
carrying a non-linear Poisson structure.

 Motivation for choosing suitable comultiplication, counit, antipode (coinverse),
and Haar weight on $(A,\Delta)$ comes from the information at the level of
Poisson--Lie groups. The proofs were given by introducing some tools like the
multiplicative unitary operator. In this way, we could argue that $(A,\Delta)$,
together with its additional structure maps, should be an example of a {\em
non-compact quantum group\/}.

 We further went on to find a ``quantum universal $R$-matrix'' type
operator related with $(A,\Delta)$ and studied its representation theory,
indicating that the ${}^*$-representations of $A$ satisfy an interesting
``quasitriangular'' type property.  See \cite{BJKp2} and \cite{BJKhj}
(More discussion on the representation theory is given in \cite{BJKppdress}.).

 However, even with these strong indications suggested by our construction
and the representation theoretic applications, we did not quite make it
clear whether $(A,\Delta)$ actually is a locally compact ($C^*$-algebraic)
quantum group.  For instance, in \cite{BJKp2}, the discussion about the
Haar weight on $(A,\Delta)$ was rather incomplete, since we restricted our
discussion to the level of a dense subalgebra of the $C^*$-algebra $A$.
Even for the (simpler) example of Rieffel's \cite{Rf5}, the full construction
of its Haar weight was not carried out.  The problem of tying together these
loose ends and establishing $(A,\Delta)$ as a locally compact quantum group
in a suitable sense was postponed to a later occasion.

 Part of the reason for the postponement was due to the fact that at the
time of writing, the question of the correct definition of a locally compact
quantum group had not yet been settled.  It was known that simply requiring
the existence of a counit and an antipode on the ``locally compact quantum
semigroup'' $(A,\Delta)$ is not enough.  Some proposals had been made, but
they were at a rather primitive stage.  Recently, the situation has improved:
A new paper by Kustermans and Vaes \cite{KuVa} appeared, in which they give
a relatively simple definition of a (reduced) $C^*$-algebraic quantum group.

 In this new definition, the existence of a left invariant (Haar) weight
and a right invariant weight plays the central role.  In particular,
they do not have to include the existence of the antipode and its polar
decomposition in their axioms.  Unlike the axiom sets of Masuda and
Nakagami \cite{MN}, or those for Kac algebras \cite{ES}, which are either
too complicated or too restrictive, these properties and others can be
proved from the defining axioms.  We are still far from achieving the goal
of formulating a set of axioms in which we do not have to invoke the existence
of Haar measure.  Considering this, it seems that the definition of Kustermans
and Vaes is the most reasonable choice at this moment.

 Now that we have an acceptable definition, we are going to return
to our example $(A,\Delta)$ and verify that $(A,\Delta)$ is indeed
a locally compact quantum group.  Our discussion will begin in Section 1
by describing the definition of Kustermans and Vaes, making precise the
notion of a ``$C^*$-algebraic locally compact quantum group''.

 In section 2, we summarize a few results about our specific example
$(A,\Delta)$.  Instead of repeating our construction carried out in
\cite{BJKp2}, we take a more economical approach of describing results
by relying less on the Poisson geometric aspects.

 In section 3,which is the main part of this paper, we describe the Haar
weight for $(A,\Delta)$ and make the notion valid in the $C^*$-algebra
setting.  Having the correct left/right invariant weights enables us to
conclude that $(A,\Delta)$ is indeed a non-compact $C^*$-algebraic
quantum group.  We are benefiting a lot from being able to work with
our specific example having a tracial weight, but many of the techniques
being used here are not necessarily type-specific, and therefore, will
be also useful in more general cases: Our discussion on the left invariance
of Haar measure is strongly motivated by the new and attractive approach
of Van Daele \cite{VDaxb}, \cite{VDoamp}.

 In Sections 4 and 5, we say a little about the antipode and the modular
function of $(A,\Delta)$.  By comparing our original definitions (motivated
by Poisson--Lie group data) with the ones suggested by the general theory,
we wish to give some additional perspectives on these maps.

 For the discussion to be complete, we need a description of the dual
counterpart to $(A,\Delta)$.  We included a very brief discussion of
$(\hat{A},\hat{\Delta})$ at the end of Section~2, and also added a short
Appendix (Section~6).  For a more careful discussion on the dual, see
\cite{BJKqhg}.  Meanwhile, we plan to pursue in our future papers the
discussion on the quantum double, as well as the research on the duality
of quantum groups in relation to the Poisson duality at the level of
their classical limit.

 As a final remark, we point out that while our original construction
of $(A,\Delta)$ was by deformation quantization of Poisson--Lie groups,
it can be also approached algebraically using the recent framework of
``twisted bicrossed products'' of Vaes and Vainerman \cite{VV}.  While
we do believe in the advantage of the more constructive approach we took
in \cite{BJKp2} motivated by Poisson--Lie groups (especially in applications
involving quantizations or representation theory, as in \cite{BJKhj}, \cite
{BJKppdress}), complementing it with the more theoretical approach presented
here will make our understanding more comprehensive.

\section{Definitions, terminologies, and conventions}

\subsection{Weights on $C^*$-algebras}
We will begin by briefly reviewing the theory of weights on $C^*$-algebras.
The purpose is to make clear the notations used in the main definition of
a $C^*$-algebraic quantum group (Definition \ref{lcqgdefn}) and in the
proofs in later sections.  For a more complete treatment and for standard
terminologies on weights, refer to \cite{Cm1}.  For instance, recall the
standard notations like ${\mathfrak N}_{\varphi}$, ${\mathfrak M}_{\varphi}$,
${{\mathfrak M}_{\varphi}}^{+}$,\dots associated to a weight $\varphi$ (on a
$C^*$-algebra $A$).  That is,
\begin{itemize}
\item ${\mathfrak N}_{\varphi}=\{a\in A:\varphi(a^*a)<\infty\}$
\item ${{\mathfrak M}_{\varphi}}^{+}=\{a\in A^+:\varphi(a)<\infty\}$
\item ${\mathfrak M}_{\varphi}={\mathfrak N}_{\varphi}^*{\mathfrak N}_{\varphi}$
\end{itemize}
The weights we will be considering are ``proper weights'': A proper weight
is a non-zero, densely defined weight on a $C^*$-algebra, which is lower
semi-continuous \cite{KuVa}.

 If we are given a (proper) weight $\varphi$ on a $C^*$-algebra, we can define
the sets ${\mathcal F}_{\varphi}$ and ${\mathcal G}_{\varphi}$ by
\begin{align}
{\mathcal F}_{\varphi}&=\{\omega\in A^*_+:\omega(x)\le\varphi(x),\forall x\in A^+\}
\notag \\
{\mathcal G}_{\varphi}&=\{\alpha\omega:\omega\in{\mathcal F}_{\varphi},\alpha
\in(0,1)\}\subseteq{\mathcal F}_{\varphi}.  \notag
\end{align}
Here $A^*$ denotes the norm dual of $A$.

These sets have been introduced by Combes, and they play a significant role
in the theory of weights.  Note that on ${\mathcal F}_{\varphi}$, one can give
a natural order inherited from $A^*_+$.  Meanwhile, ${\mathcal G}_{\varphi}$ is
a directed subset of ${\mathcal F}_{\varphi}$.  That is, for every $\omega_1,
\omega_2\in{\mathcal G}_{\varphi}$, there exists an element $\omega\in
{\mathcal G}_{\varphi}$ such that $\omega_1,\omega_2\le\omega$.  Because
of this, ${\mathcal G}_{\varphi}$ is often used as an index set (of a net). 
For a proper weight $\varphi$, we would have: $\varphi(x)=\lim\bigl(\omega(x)
\bigr)_{\omega\in{\mathcal G}_{\varphi}}$, for $x\in A^+$.

 By standard theory, for a weight $\varphi$ on a $C^*$-algebra $A$, one can
associate to it a ``GNS-construction'' $({\mathcal H}_{\varphi},{\pi}_{\varphi},
{\Lambda}_{\varphi})$.  Here, ${\mathcal H}_{\varphi}$ is a Hilbert space,
${\Lambda}_{\varphi}:{\mathfrak N}_{\varphi}\to{\mathcal H}_{\varphi}$ is
a linear map such that ${\Lambda}_{\varphi}({\mathfrak N}_{\varphi})$ is
dense in ${\mathcal H}_{\varphi}$ and $\bigl\langle{\Lambda}_{\varphi}
(a),{\Lambda}_{\varphi}(b)\bigr\rangle=\varphi(b^*a)$ for $a,b\in
{\mathfrak N}_{\varphi}$, and $\pi_{\varphi}$ is a representation of $A$
on ${\mathcal H}_{\varphi}$ defined by $\pi_{\varphi}(a)\Lambda_{\varphi}
(b)=\Lambda_{\varphi}(ab)$ for $a\in A$, $b\in{\mathfrak N}_{\varphi}$.
The GNS-construction is unique up to a unitary
transformation.

 If $\varphi$ is proper, then ${\mathfrak N}_{\varphi}$ is dense in $A$
and ${\Lambda}_{\varphi}:{\mathfrak N}_{\varphi}\to{\mathcal H}_{\varphi}$
is a closed map.  Also $\pi_{\varphi}:A\to{\mathcal B}({\mathcal H}_{\varphi})$
is a non-degenerate ${}^*$-homomorphism.  It is not difficult to show that
$\varphi$ has a natural extension to a weight on $M(A)$, which we will still
denote by $\varphi$.  Meanwhile, since we can define for every $\omega\in
{\mathcal G}_{\varphi}$ a unique element $\tilde{\omega}\in\pi_{\varphi}(A)''
{}_*$ such that $\tilde{\omega}\circ\pi_{\varphi}=\omega$, we can define
a weight $\tilde{\varphi}$ on the von Neumann algebra $\pi_{\varphi}(A)''$
in the following way: $\tilde{\varphi}(x)=\lim\bigl(\tilde{\omega}(x)
\bigr)_{\omega\in{\mathcal G}_{\varphi}}$, for $x\in\bigl(\pi_{\varphi}(A)''
\bigr)^+$.  Then, by standard terminology \cite{St}, $\tilde{\varphi}$
is a ``normal'', ``semi-finite'' weight on the von Neumann algebra
$\pi_{\varphi}(A)''$.

 Motivated by the properties of normal, semi-finite weights on von
Neumann algebras, and to give somewhat of a control over the
non-commutativity of $A$, one introduces the notion of ``KMS weights''
\cite{Ku1}.  The notion as defined below is slightly different from
(but equivalent to) the original one given by Combes in \cite{Cm2}.

\begin{defn}\label{KMS}
A proper weight $\varphi$ is called a ``KMS weight'' if there exists a
norm-continuous one-parameter group of automorphisms $\sigma$
of $A$ such that
\begin{enumerate}
\item $\varphi\circ\sigma_t=\varphi$, for all $t\in{\mathbb{R}}$.
\item $\varphi(a^*a)=\varphi\bigl(\sigma_{i/2}(a)\sigma_{i/2}(a)^*\bigr)$,
for all $a\in D(\sigma_{i/2})$.
\end{enumerate}
Here $\sigma_{i/2}$ is the analytic extension of the one-parameter group
$\sigma_t$ to $\frac{i}{2}$.
\end{defn}
The one-parameter group $\sigma$ is called the ``modular automorphism group''
for $\varphi$.  It is uniquely determined when $\varphi$ is faithful.
Meanwhile, a proper weight $\varphi$ is said to be ``approximately KMS''
if the associated (normal, semi-finite) weight $\tilde{\varphi}$ is
faithful. A KMS weight is approximately KMS.  For more discussion
on these  classes of weights, including the relationship between the
conditions above and the usual KMS condition, see \cite{Ku1}.  Finally,
note that in the special case when $\varphi$ is a trace (i.\,e. $\varphi
(a^*a)=\varphi(aa^*)$, for $a\in{\mathfrak N}_{\varphi}$), it is clear
that $\varphi$ is KMS.  The modular automorphism group will be trivial
($\equiv\operatorname{Id}$).

\subsection{Definition of a locally compact quantum group}
Let $A$ be a $C^*$-algebra.  Suppose $\Delta:A\to M(A\otimes A)$ is
a non-degenerate ${}^*$-homomorphism (Later, $\Delta$ will be given
certain conditions for it to become a comultiplication.).  A proper weight
$\varphi$ on $(A,\Delta)$ will be called {\em left invariant\/}, if
\begin{equation}\label{(left)}
\varphi\bigl((\omega\otimes\operatorname{id})(\Delta a)\bigr)
=\omega(1)\varphi(a),
\end{equation}
for all $a\in{{\mathfrak M}_{\varphi}}^+$ and $\omega\in A^*_+$.  Similarly,
$\varphi$ is called {\em right invariant\/}, if
\begin{equation}
\varphi\bigl((\operatorname{id}\otimes\omega)(\Delta a)\bigr)
=\omega(1)\varphi(a).
\end{equation}
By $\omega(1)$, we mean $\|\omega\|$.  Note here that we used the
extensions of $\varphi$ to $M(A)$ in the equations, since we only know
that $(\omega\otimes\operatorname{id})(\Delta a)\in M(A)^+$.  In the
definition of locally compact quantum groups (to be given below), the
``slices'' of $\Delta a$ will be assumed to be contained in $A$.

 In the definitions above, the left [respectively, right] invariance
condition requires the formula \eqref{(left)} to hold only for $a\in
{{\mathfrak M}_{\varphi}}^+$.  It is a very weak form of left invariance.
In the case of locally compact quantum groups, the result can be extended
and a much stronger left invariance condition can be proved from it.
The proof is non-trivial.  It was one of the important contributions
made by Kustermans and Vaes.

 Next, let us state the definition of a locally compact ($C^*$-algebraic)
quantum group given by Kustermans and Vaes \cite{KuVa}.  In the
definition, $[X]$ denotes the closed linear span of $X$.
\begin{defn}\label{lcqgdefn}
Consider a $C^*$-algebra $A$ and a non-degenerate
${}^*$-homomorphism $\Delta:A\to M(A\otimes A)$ such that
\begin{enumerate}
\item $(\Delta\otimes\operatorname{id})\Delta=(\operatorname{id}
\otimes\Delta)\Delta$
\item $\bigl[\bigl\{(\omega\otimes\operatorname{id})(\Delta a):
\omega\in A^*,a\in A\bigr\}\bigr]=A$
\item $\bigl[\bigl\{(\operatorname{id}\otimes\omega)(\Delta a):
\omega\in A^*,a\in A\bigr\}\bigr]=A$
\end{enumerate}
Moreover, assume that there exist weights $\varphi$ and $\psi$ such that
\begin{itemize}
\item $\varphi$ is a faithful, left invariant approximate KMS weight on
$(A,\Delta)$.
\item $\psi$ is a right invariant approximate KMS weight on $(A,\Delta)$.
\end{itemize}
Then we say that $(A,\Delta)$ is a {\em (reduced) $C^*$-algebraic
quantum group\/}.
\end{defn}

 First condition is the ``coassociativity'' condition for the
``comultiplication'' $\Delta$.  By the non-degeneracy, it can
be naturally extended to $M(A)$ [we can also extend $(\Delta\otimes
\operatorname{id})$ and $(\operatorname{id}\otimes\Delta)$], thereby
making the expression valid.  The two density conditions more or less
correspond to the cancellation property in the case of ordinary groups,
although they are somewhat weaker.  The last axiom corresponds to
the existence of Haar measure (The weights $\varphi$ and $\psi$
actually turn out to be faithful KMS weights.).  For more on this
definition (e.\,g. discussions on how one can build other structure
maps like the antipode), see \cite{KuVa}.

\section{The Hopf $C^*$-algebra $(A,\Delta)$}

 Our main object of study is the Hopf $C^*$-algebra $(A,\Delta)$
constructed in \cite{BJKp2}.  As a $C^*$-algebra, $A$ is isomorphic to
a twisted group $C^*$-algebra $C^*\bigl(H/Z,C_0(\Gg/\Gq),\sigma\bigr)$,
where $H$ is the $(2n+1)$-dimensional Heisenberg Lie group and $Z$ is
the center of $H$.  Whereas, $\Gg=\Gh^*$ is the dual space of the Lie
algebra $\Gh$ of $H$ and $\Gq=\Gz^{\bot}$, for $\Gz\subseteq\Gh$
corresponding to $Z$.  Since $H$ is a nilpotent Lie group, $H\cong\Gh$
and $Z\cong\Gz$, as vector spaces.  We denoted by $\sigma$ (not to be
confused with the modular automorphism group) the twisting cocycle
for the group $H/Z$.  As constructed in \cite{BJKp2}, $\sigma$ is
a continuous field of cocycles $\Gg/\Gq\ni r\mapsto\sigma^r$, where
\begin{equation}\label{(sigma)}
\sigma^r\bigl((x,y),(x',y')\bigr)=\bar{e}\bigl[\eta_{\lambda}(r)
\beta(x,y')\bigr].
\end{equation}
Following the notation of the previous paper, we used: $\bar{e}(t)
=e^{(-2\pi i)t}$ and $\eta_{\lambda}(r)=\frac{e^{2\lambda r}-1}
{2\lambda}$, where $\lambda$ is a fixed real constant.  We denote
by $\beta(\ ,\ )$ the inner product.  The elements $(x,y)$, $(x',y')$
are group elements in $H/Z$.

 In \cite{BJKp2}, we showed that the $C^*$-algebra $A$ is a deformation
quantization (in Rieffel's ``strict'' sense \cite{Rf1}, \cite {Rf4}) of
$C_0(G)$, where $G$ is a certain solvable Lie group which is the {\em
dual Poisson--Lie group\/} of $H$.  The number $\lambda$ mentioned above
determines the group structure of $G$ (When $\lambda=0$, the group $G$
becomes abelian, which is not very interesting.).  See Definition 1.6 of
\cite{BJKp2} for the precise definition of $G$.  For convenience, we fixed
the deformation parameter as $\hbar=1$. This is the reason why we do not
see it in the definition of $A$. If we wish to illustrate the deformation
process, we may just replace $\beta$ by $\hbar\beta$, and let $\hbar\to0$.
When $\hbar=0$ (i.\,e. classical limit), we have $\sigma\equiv1$, and hence,
$A_{\hbar=0}\cong C_0(G)$.  Throughout this paper, we will just work with
$A=A_{\hbar=1}$.

 Let us be a little more specific and recall some of the
notations and results obtained in \cite{BJKp2}, while referring
the reader to that paper for more details on the construction
of our main example $(A,\Delta)$.

 We first introduce the subspace ${\mathcal A}$, which is a dense
subspace of $A$ consisting of the functions in $S_{3c}(H/Z\times
\Gg/\Gq)$, the space of Schwartz functions in the $(x,y,r)$
variables having compact support in the $r(\in\Gg/\Gq)$ variable.
On ${\mathcal A}$, we define the (twisted) multiplication and the
(twisted) involution as follows:
\begin{equation}\label{(multiplication)}
(f\times g)(x,y,r)=\int f(\tilde{x},\tilde{y},r)g(x-\tilde{x},
y-\tilde{y},r)\bar{e}\bigl[\eta_{\lambda}(r)\beta(\tilde{x},
y-\tilde{y})\bigr]\,d\tilde{x}d\tilde{y},
\end{equation}
and
\begin{equation}\label{(involution)}
f^*(x,y,r)=\overline{f(-x,-y,r)}\bar{e}\bigl[\eta_{\lambda}(r)
\beta(x,y)\bigr].
\end{equation}
It is not difficult to see that ${\mathcal A}=S_{3c}(H/Z\times\Gg/
\Gq)$ is closed under the multiplication \eqref{(multiplication)}
and the involution \eqref{(involution)}.  Here, we observe the
role being played by the twisting cocycle $\sigma$ defined in
\eqref{(sigma)}.

 Elements of ${\mathcal A}$ are viewed as operators on the Hilbert
space ${\mathcal H}=L^2(H/Z\times\Gg/\Gq)$, via the ``regular
representation'', $L$, defined by
\begin{equation}\label{(representation)}
(L_f\xi)(x,y,r)=\int f(\tilde{x},\tilde{y},r)\xi(x-\tilde{x},
y-\tilde{y},r)\bar{e}\bigl[\eta_{\lambda}(r)\beta(\tilde{x},
y-\tilde{y})\bigr]\,d\tilde{x}d\tilde{y}.
\end{equation}
For $f\in{\mathcal A}$, define its norm by $\|f\|=\|L_f\|$.
Then $({\mathcal A},\times,{}^*,\|\ \|)$ as above is a
pre-$C^*$-algebra, whose completion is the $C^*$-algebra $A
\cong C^*\bigl(H/Z,C_0(\Gg/\Gq),\sigma\bigr)$.

\begin{rem}
To be more precise, the completion of ${\mathcal A}$ with respect to
the norm given by the regular representation, $L$, should be isomorphic
to the ``reduced'' twisted group $C^*$-algebra $C^*_r\bigl(H/Z,C_0
(\Gg/\Gq),\sigma\bigr)$.  But by using a result of Packer and Raeburn
\cite{PR}, it is rather easy to see that the amenability condition
holds in our case, thereby obtaining the isomorphism with the ``full''
$C^*$-algebra as above.  Meanwhile, we should point out that our
definition of ${\mathcal A}$ is slightly different from that of
\cite{BJKp2}: There, ${\mathcal A}$ is a subspace of $C_0(G)$, while
at present we view it as functions contained in $C_0(H/Z\times\Gg/\Gq)$,
in the $(x,y,r)$ variables.  Nevertheless, they can be regarded as the
same since we consider the functions in ${\mathcal A}$ as operators
contained in our $C^*$-algebra $A$.  The identification of the function
spaces is given by the (partial) Fourier transform.
\end{rem}

 The $C^*$-algebra $A$ becomes a Hopf $C^*$-algebra, together with its
{\em comultiplication\/} $\Delta$.  In the following proposition, we chose
to describe the comultiplication in terms of a certain ``multiplicative
unitary operator'' $U_A\in{\mathcal B}({\mathcal H}\otimes{\mathcal H})$.
See \cite{BJKp2} for a discussion on the construction of $U_A$.

\begin{prop}\label{comultiplication}
\begin{enumerate}
\item Let $U_A$ be the operator on ${\mathcal H}\otimes{\mathcal H}$
defined by
\begin{align}
U_A\xi(x,y,r,x',y',r')&=(e^{-\lambda r'})^n\bar{e}\bigl[\eta_{\lambda}
(r')\beta(e^{-\lambda r'}x,y'-e^{-\lambda r'}y)\bigr] \notag \\
&\quad\xi(e^{-\lambda r'}x,e^{-\lambda r'}y,r+r',
x'-e^{-\lambda r'}x,y'-e^{-\lambda r'}y,r').  \notag
\end{align}
Then $U_A$ is a unitary operator, and is multiplicative.  That is,
$$
U_{12}U_{13}U_{23}=U_{23}U_{12}.
$$
\item For $f\in{\mathcal A}$, define $\Delta f$ by
$$
\Delta f=U_A(f\otimes1){U_A}^*,
$$
where $f$ and $\Delta f$ are understood as operators $L_f$ and
$(L\otimes L)_{\Delta f}$.  Then $\Delta$ can be extended to a
non-degenerate $C^*$-homomorphism $\Delta:A\to M(A\otimes A)$
satisfying the coassociativity condition:
$$
(\Delta\otimes\operatorname{id})(\Delta f)=(\operatorname{id}
\otimes\Delta)(\Delta f).
$$
\end{enumerate}
\end{prop}

\begin{proof}
See Proposition 3.1 and Theorem 3.2 of \cite{BJKp2}, together
with the Remark 3.3 following them.
\end{proof}

 There is a useful characterization of the $C^*$-algebra $A$, via
the multiplicative unitary operator $U_A$.  The following result
is suggested by the general theory on multiplicative unitaries
by Baaj and Skandalis \cite{BS}.

\begin{prop}\label{A(U)}
Let $U_A$ be as above.  Consider the subspace ${\mathcal A}(U_A)$
of ${\mathcal B}({\mathcal H})$ defined below:
$$
{\mathcal A}(U_A)=\bigl\{(\omega\otimes\operatorname{id})(U_A):
\omega\in{\mathcal B}({\mathcal H})_*\bigr\}.
$$
By standard theory, ${\mathcal A}(U_A)$ is a subalgebra of the
operator algebra ${\mathcal B}({\mathcal H})$, and the subspace
${\mathcal A}(U_A){\mathcal H}$ forms a total set in ${\mathcal H}$.

We can show that the norm-closure in ${\mathcal B}({\mathcal H})$
of the algebra ${\mathcal A}(U_A)$ is exactly the $C^*$-algebra $A$
we are studying.  That is,
$$
A=\overline{\bigl\{(\omega\otimes\operatorname{id})(U_A):\omega\in
{\mathcal B}({\mathcal H})_*\bigr\}}^{\|\ \|}.
$$
\end{prop}

\begin{proof}
The definition and the properties of ${\mathcal A}(U_A)$ can be
found in \cite{BS}.  We only need to verify the last statement.
We will work with the standard notation $\omega_{\xi,\eta}$, where
$\xi,\eta\in{\mathcal H}$.  It is defined by $\omega_{\xi,\eta}
(a)=\langle a\xi,\eta\rangle$, and it is well known that linear
combinations of the $\omega_{\xi,\eta}$ are (norm) dense in
${\mathcal B}({\mathcal H})_*$.

So consider $(\omega_{\xi,\eta}\otimes\operatorname{id})(U_A)\in
{\mathcal B}({\mathcal H})$.  We may further assume that $\xi$
and $\eta$ are continuous functions having compact support.
Let $\zeta\in{\mathcal H}$.  Then, by using change of variables,
we have:
\begin{align}
&(\omega_{\xi,\eta}\otimes\operatorname{id})(U_A)\zeta(x,y,r)
\notag \\
&=\int\bigl(U_A(\xi\otimes\zeta)\bigr)(\tilde{x},\tilde{y},\tilde{r};
x,y,r)\overline{\eta(\tilde{x},\tilde{y},\tilde{r})}\,d\tilde{x}
d\tilde{y}d\tilde{r}  \notag \\
&=\int f(\tilde{x},\tilde{y},r)\zeta(x-\tilde{x},y-\tilde{y},r)
\bar{e}\bigl[\eta_{\lambda}(r)\beta(\tilde{x},y-\tilde{y})
\bigr]\,d\tilde{x}d\tilde{y},  \notag
\end{align}
where
$$
f(\tilde{x},\tilde{y},r)=\int\xi(\tilde{x},\tilde{y},\tilde{r}+r)
(e^{\lambda r})^n\overline{\eta(e^{\lambda r}\tilde{x},e^{\lambda r}
\tilde{y},\tilde{r})}\,d\tilde{r}.
$$
Since $\xi$ and $\eta$ are $L^2$-functions, the integral (and thus
$f$) is well defined.  Actually, since $f$ is essentially defined
as a convolution product (in $r$) of two continuous functions having
compact support, $f$ will be also continuous with compact support.
This means that
$$
(\omega_{\xi,\eta}\otimes\operatorname{id})(U_A)=L_f\in A.
$$
Meanwhile, since the choice of $\xi$ and $\eta$ is arbitrary, we
can see that the collection of the $f$ will form a total set in
the space of continuous functions in the $(x,y,r)$ variables having
compact support.  It follows from these two conclusions that
$$
\overline{\bigl\{(\omega\otimes\operatorname{id})(U_A):\omega\in
{\mathcal B}({\mathcal H})_*\bigr\}}^{\|\ \|}=A.
$$
\end{proof}

 Meanwhile, from the proof of Theorem 3.2 of \cite{BJKp2}, we also have
the following result.  These are not same as the density conditions of
Definition \ref{lcqgdefn}, but are actually stronger: This is rather
well known and can be seen easily by applying linear functionals (Use the
fact that any $\omega\in A^*$ has the form $\omega'(\cdot\,b)$, with $\omega'
\in A^*$ and $b\in A$.).

\begin{prop}\label{density}
We have:
\begin{enumerate}
\item $\Delta (A)(1\otimes A)$ is dense in $A\otimes A$.
\item $\Delta (A)(A\otimes 1)$ is dense in $A\otimes A$.
\end{enumerate}
\end{prop}

\begin{proof}
In the proof of the non-degeneracy of $\Delta$ in Theorem 3.2
and Remark 3.3 of \cite{BJKp2}, we showed that the $(\Delta f)
(1\otimes g)$'s (for $f,g\in{\mathcal A}$) form a total set in
the space $S_{3c}(H/Z\times\Gg/\Gq\times H/Z\times\Gg/\Gq)$, which
is in turn shown to be dense in $A\otimes A$: Under the natural
injection from $S_{3c}(H/Z\times\Gg/\Gq\times H/Z\times\Gg/\Gq)$
into ${\mathcal B}({\mathcal H}\otimes{\mathcal H})$, the algebraic
tensor product ${\mathcal A}\odot{\mathcal A}$ is sent into a
dense subset of the algebraic tensor product $A\odot A$.  Since
elements in $S_{3c}(H/Z\times\Gg/\Gq\times H/Z\times\Gg/\Gq)$ can
be approximated (in the $L^1$-norm) by elements of ${\mathcal A}
\odot{\mathcal A}$, we see that $S_{3c}(H/Z\times\Gg/\Gq\times H/Z
\times\Gg/\Gq)$ is mapped into a dense subset (in the $C^*$ norm)
of $A\otimes A$.  Thus it follows that $\Delta (A)(1\otimes A)$
is dense in $A\otimes A$.  The second statement can be shown in
exactly the same way.
\end{proof}

 Turning our attention to the other structures on $(A,\Delta)$,
we point out that by viewing $A$ as a ``quantum $C_0(G)$'', we
can construct its {\em counit\/}, $\varepsilon$, and {\em antipode\/},
$S$.  These are described in the following proposition.

\begin{prop}\label{antipode}
\begin{enumerate}
\item For $f\in{\mathcal A}$, define $\varepsilon:{\mathcal A}\to
\mathbb{C}$ by
$$
\varepsilon(f)=\int f(x,y,0)\,dxdy.
$$
Then $\varepsilon$ can be extended to a $C^*$-homomorphism from $A$
to $\mathbb{C}$ satisfying the condition: $(\operatorname{id}
\otimes\varepsilon)\Delta=(\varepsilon\otimes\operatorname{id})\Delta
=\operatorname{id}$.
\item Consider a map $S:{\mathcal A}\to{\mathcal A}$ defined by
$$
\bigl(S(f)\bigr)(x,y,r)=(e^{2\lambda r})^{n}\bar{e}\bigl[\eta_{\lambda}(r)
\beta(x,y)\bigr]f(-e^{\lambda r}x,-e^{\lambda r}y,-r).
$$
Then $S$ can be extended to an anti-automorphism $S:A\to A$, satisfying:
$S\bigl(S(a)^*\bigr)^*=a$ and $(S\otimes S)(\Delta a)=\chi\bigl(\Delta(S(a))
\bigr)$, where $\chi$ denotes the flip.  Actually, we have:
$S^2=\operatorname{Id}$.
\end{enumerate}
\end{prop}

\begin{proof}
See Theorem 4.1 and Proposition 4.3 of \cite{BJKp2}.  We had to use partial
Fourier transform to convert these results into the level of functions in
the $(x,y,r)$ variables.  We also mention here that $S$ is defined by $S(a)=
\hat{J}a^*\hat{J}$, where $\hat{J}$ is an involutive operator on ${\mathcal H}$
defined by
$$
\hat{J}\xi(x,y,r)=(e^{\lambda r})^n\overline{\xi(e^{\lambda r}x,
e^{\lambda r}y,-r)}.
$$
Since $\hat{J}$ is an anti-unitary involutive operator, it is easy to see
that $S$ is an anti-automorphism such that $S^2=\operatorname{Id}$.
\end{proof}

\begin{rem}
The notation for the operator $\hat{J}$ is motivated by the modular theory
and by \cite{MN}.  Meanwhile, since the square of the antipode is the identity,
our example is essentially the {\em Kac $C^*$-algebra\/} (Compare with \cite{EV},
although our example is actually non-unimodular, unlike in that paper.  See also
\cite{VV}, of which our example is a kind of a special case.).  We also note
that in \cite{BJKp2}, we used $\kappa$ to denote the antipode, while we use $S$
here.  This is done so that we can match our notation with the preferred notation
of Kustermans and Vaes \cite{KuVa}.
\end{rem}

 In general, the counit may as well be unbounded.  So Proposition~\ref
{antipode} implies that what we have is a more restrictive ``bounded
counit''.  Having $S$ bounded is also a bonus.  Even so, the result
of the proposition is not enough to legitimately call $S$ an antipode.
To give some support for our choice, we also showed the following,
albeit only at the level of the function space $\mathcal A$.  See
section 4 of \cite{BJKp2}.

\begin{prop}\label{antipodem}
For $f\in{\mathcal A}$, we have:
$$
m\bigl((\operatorname{id}\otimes S)(\Delta f)\bigr)
=m\bigl((S\otimes\operatorname{id})(\Delta f)\bigr)=\varepsilon(f)1,
$$
where $m:{\mathcal A}\otimes{\mathcal A}\to{\mathcal A}$ is the
multiplication.
\end{prop}

 This is the required condition for the antipode in the purely algebraic
setting of Hopf algebra theory \cite{Mo}.  In this sense, the proposition
gives us a modest justification for our choice of $S$.  However, in the
operator algebra setting, this is not the correct way of approach.  One
of the serious obstacles is that the multiplication $m$ is in general not
continuous for the operator norm, thereby giving us trouble extending $m$
to $A\otimes A$ or $M(A\otimes A)$.

 Because of this and other reasons (including the obstacles due to
possible unboundedness of $\varepsilon$ and $S$), one has to develop
a new approach.  Motivated by the theory of Kac algebras \cite{ES},
operator algebraists have been treating the antipode together with
the notion of the Haar weight.  This is also the approach chosen by
Masuda, Nakagami \cite{MN} and by Kustermans, Vaes \cite{KuVa}.  As we
mentioned earlier in this paper, any rigorous discussion about locally
compact quantum groups should be built around the notion of Haar
weights.  In the next section, we will exclusively discuss the Haar
weight for our $(A,\Delta)$, and establish that $(A,\Delta)$ is indeed
a ``$C^*$-algebraic locally compact quantum group''.  We will come back
to the discussion of the antipode in section 4.

 Before wrapping up this section, let us mention the dual object for our
$(A,\Delta)$, which would be the ``dual locally compact quantum group''
$(\hat{A},\hat{\Delta})$.  Our discussion here is kept to a minimum.
More careful discussion on $(\hat{A},\hat{\Delta})$ is presented in a
separate paper \cite{BJKqhg}.  Meanwhile, see Appendix (Section 6) for
a somewhat different characterization of the dual object.

\begin{prop}\label{dualqg}
\begin{enumerate}
\item Let $U_A$ be as above.  Consider the subspace $\hat{\mathcal A}
(U_A)$ of ${\mathcal B}({\mathcal H})$ defined by
$$
\hat{\mathcal A}(U_A)=\bigl\{(\operatorname{id}\otimes\omega)(U_A):
\omega\in{\mathcal B}({\mathcal H})_*\bigr\}.
$$
Then $\hat{\mathcal A}(U_A)$ is a subalgebra of the operator algebra
${\mathcal B}({\mathcal H})$, and the subspace $\hat{\mathcal A}(U_A)
{\mathcal H}$ forms a total set in ${\mathcal H}$.
\item The norm-closure in ${\mathcal B}({\mathcal H})$ of the algebra $\hat
{\mathcal A}(U_A)$ is the $C^*$-algebra $\hat{A}$:
$$
\hat{A}=\overline{\bigl\{(\operatorname{id}\otimes\omega)(U_A):\omega\in
{\mathcal B}({\mathcal H})_*\bigr\}}^{\|\ \|}.
$$
\item For $b\in\hat{\mathcal A}(U_A)$, define $\hat{\Delta}b$ by
$\hat{\Delta}b={U_A}^*(1\otimes b)U_A$.  Then $\hat{\Delta}$ can be
extended to the comultiplication $\hat{\Delta}:\hat{A}\to M(\hat{A}
\otimes\hat{A})$.
\item The duality exists at the level of ${\mathcal A}(U_A)$ and $\hat
{\mathcal A}(U_A)$, by the following formula:
$$
\bigl\langle L(\omega),\rho({\omega}')\bigr\rangle=(\omega\otimes
{\omega}')(U_A)=\omega\bigl(\rho({\omega}')\bigr)={\omega}'\bigl(
L(\omega)\bigr),
$$
where $L(\omega)=(\omega\otimes\operatorname{id})(U_A)\in{\mathcal A}(U_A)$
and $\rho({\omega}')=(\operatorname{id}\otimes{\omega}')(U_A)\in\hat{\mathcal A}
(U_A)$.
\end{enumerate}
\end{prop}

\begin{proof}
Since our multiplicative unitary operator $U_A$ is known to be regular,
we can just follow the standard theory of multiplicative unitary operators
\cite{BS}.
\end{proof}

 At this moment, $(\hat{A},\hat{\Delta})$ is just a quantum semigroup,
having only the comultiplication.  However, once we establish in section 3
the proof that our $(A,\Delta)$ is indeed a ($C^*$-algebraic) locally compact
quantum group, we can apply the general theory \cite{KuVa}, and show that
$(\hat{A},\hat{\Delta})$ is also a locally compact quantum group.  Meanwhile,
we can give a more specific description of the $C^*$-algebra $\hat{A}$, as
presented below.  

\begin{prop}\label{rhof}
Let $\hat{\mathcal A}$ be the space of Schwartz functions in the $(x,y,r)$
variables having compact support in the $r$ variable.  For $f\in\hat{\mathcal A}$,
define the operator ${\rho}_f\in{\mathcal B}({\mathcal H})$ by
\begin{equation}\label{(rhof)}
({\rho}_f\zeta)(x,y,r)=\int(e^{\lambda\tilde{r}})^n f(x,y,\tilde{r})
\zeta(e^{\lambda\tilde{r}}x,e^{\lambda\tilde{r}}y,r-\tilde{r})\,d\tilde{r}.
\end{equation}
Then the $C^*$-algebra $\hat{A}$ is generated by the operators ${\rho}_f$.
\end{prop}

\begin{proof}
As in the proof of Proposition \ref{A(U)}, consider the operators
$(\operatorname{id}\otimes\omega_{\xi,\eta})(U_A)$ in $\hat{\mathcal A}
(U_A)$.  Without loss of generality, we can assume that $\xi$ and
$\eta$ are continuous functions having compact support.  Let $\zeta
\in{\mathcal H}$.  Then we have:
\begin{align}
&(\operatorname{id}\otimes\omega_{\xi,\eta})(U_A)\zeta(x,y,r)
\notag \\
&=\int\bigl(U_A(\zeta\otimes\xi)\bigr)(x,y,r;\tilde{x},\tilde{y},\tilde{r})
\overline{\eta(\tilde{x},\tilde{y},\tilde{r})}\,d\tilde{x}d\tilde{y}d\tilde{r}
\notag \\
&=\int(e^{\lambda\tilde{r}})^n f(x,y,\tilde{r})\zeta(e^{\lambda\tilde{r}}x,
e^{\lambda\tilde{r}}y,r-\tilde{r})\,d\tilde{r}, \notag
\end{align}
where
$$
f(x,y,\tilde{r})=\int\bar{e}\bigl[\eta_{\lambda}(\tilde{r})\beta
(x,y-e^{-\lambda\tilde{r}}\tilde{y})\bigr]\xi(\tilde{x}-e^{\lambda
\tilde{r}}x,\tilde{y}-e^{\lambda\tilde{r}}y,-\tilde{r})\overline
{\eta(\tilde{x},\tilde{y},-\tilde{r})}\,d\tilde{x}d\tilde{y}.
$$
We see that $(\operatorname{id}\otimes\omega_{\xi,\eta})(U_A)={\rho}_f$,
and $f$ is continuous with compact support.  By the same argument we
used in the proof of Proposition \ref{A(U)}, we conclude that the
$C^*$-algebra generated by the operators ${\rho}_f$, $f\in\hat{\mathcal A}$,
coincides with the $C^*$-algebra $\hat{A}$.
\end{proof}

 The above characterization of $\hat{A}$ is useful when we wish to
regard $\hat{A}$ as a deformation quantization of a Poisson--Lie group.
Actually, by using partial Fourier transform, we can show without
difficulty that $\hat{A}\cong\rho\bigl(C^*(G)\bigr)$, where $\rho$
is the right regular representation of $C^*(G)$.  In a future paper,
we will have an occasion to discuss the duality between $(A,\Delta)$
and $(\hat{A},\hat{\Delta})$ in relation to the Poisson--Lie group
duality between $G$ and $H$.

\section{Haar weight}

 We have been arguing that $(A,\Delta)$ {\em is\/} a ``quantized
$C_0(G)$''.  This viewpoint has been helpful in our construction of
its comultiplication $\Delta$, counit $\varepsilon$, and antipode $S$.

 To discuss the (left invariant) Haar weight on $(A,\Delta)$, we
pull this viewpoint once more.  Recall that the group structure on
$G$ has been specifically chosen (Definition 1.6 of \cite{BJKp2})
so that the Lebesgue measure on $G$ becomes its left invariant Haar
measure.  This suggests us to build the Haar weight on $(A,\Delta)$
from the Lebesgue measure on $G$.  At the level of functions in
$\mathcal A$, this suggestion takes the following form:

\begin{defn}\label{Haar}
On ${\mathcal A}$, define a linear functional $\varphi$ by
$$
\varphi(f)=\int f(0,0,r)\,dr.
$$
\end{defn}

 In section 5 of \cite{BJKp2}, we obtained some results (including
the ``left invariance'' property) indicating that our choice of
$\varphi$ is a correct one.  However, the discussion was limited to
the level of functions in $\mathcal A$, and thus not very satisfactory.

 Jumping up from the function level to the operator level can be quite
technical, and it is not necessarily an easy task (For example,
see \cite{Bj}, \cite{VDaxb}.).  Whereas, if one wants to rigorously
formulate the construction of a locally compact quantum group in the
operator algebra setting, this step of ``jumping up'' (extending
$\varphi$ to a weight) is very crucial.

 Fortunately in our case, the discussion will be much simpler than some
of the difficult examples, since we can show that $\varphi$ is tracial.
Note the following:

\begin{prop}\label{phitrace}
Let $\varphi$ be defined on $\mathcal A$ as in Definition \ref{Haar}.
Then for $f\in{\mathcal A}$, we have:
$$
\varphi(f^*\times f)=\varphi(f\times f^*)=\|f\|_2^2,
$$
where $f^*$ is the $C^*$-involution of $f$, as given in \eqref{(involution)}.
\end{prop}

\begin{proof}
By equations \eqref{(multiplication)} and \eqref{(involution)}, we have:
\begin{align}
&(f^*\times f)(x,y,r)=\int f^*(\tilde{x},\tilde{y},r)f(x-\tilde{x},
y-\tilde{y},r)\bar{e}\bigl[\eta_{\lambda}(r)\beta(\tilde{x},y-\tilde{y})
\bigr]\,d\tilde{x}d\tilde{y}   \notag \\
&=\int\overline{f(-\tilde{x},-\tilde{y},r)}\bar{e}\bigl[\eta_{\lambda}(r)
\beta(\tilde{x},\tilde{y})\bigr]f(x-\tilde{x},y-\tilde{y},r)\bar{e}\bigl[
\eta_{\lambda}(r)\beta(\tilde{x},y-\tilde{y})\bigr]\,d\tilde{x}d\tilde{y}
\notag \\
&=\int\overline{f(-\tilde{x},-\tilde{y},r)}f(x-\tilde{x},y-\tilde{y},r)
\bar{e}\bigl[\eta_{\lambda}(r)\beta(\tilde{x},y)\bigr]\,d\tilde{x}d\tilde{y}.
\notag
\end{align}
It follows that:
\begin{align}
\varphi(f^*\times f)&=\int (f^*\times f)(0,0,r)\,dr
=\int\overline{f(-\tilde{x},-\tilde{y},r)}f(-\tilde{x},-\tilde{y},r)
\,d\tilde{x}d\tilde{y}dr  \notag \\
&=\int\overline{f(\tilde{x},\tilde{y},r)}f(\tilde{x},\tilde{y},r)
\,d\tilde{x}d\tilde{y}dr=\|f\|_2^2. \notag 
\end{align}
The identity $\varphi(f\times f^*)=\|f\|_2^2$ can be proved similarly.
\end{proof}

\begin{cor}
By Proposition \ref{phitrace}, we see that $\varphi$ is a faithful,
positive linear functional which is a trace.
\end{cor}

 Now, let us begin the discussion of constructing a weight on $(A,
\Delta)$ extending $\varphi$.  As a first step, let us consider the
associated GNS construction for $\varphi$.  We can see below that the
``regular representation'' $L$ on ${\mathcal H}$ we defined earlier
is the GNS representation for $\varphi$.

\begin{prop}
Consider the Hilbert space ${\mathcal H}=L^2(H/Z\times\Gg/\Gq)$,
and let $\Lambda:{\mathcal A}\hookrightarrow{\mathcal H}$ be the
inclusion map.  Then for $f,g\in{\mathcal A}$, we have:
$$
\bigl\langle\Lambda(f),\Lambda(g)\bigr\rangle
=\varphi(g^*\times f).
$$
Here $\langle\ ,\ \rangle$ is the inner product on $\mathcal H$, conjugate
in the second place.  Meanwhile, left multiplication gives a non-degenerate
${}^*$-representation, $\pi_{\varphi}:{\mathcal A}\to{\mathcal B}({\mathcal H})$,
which coincides with the ``regular representation'' $L$.

By essential uniqueness of the GNS construction, we conclude that $({\mathcal H},
\Lambda,\pi_{\varphi})$ is the GNS triple associated with $\varphi$.
\end{prop}

\begin{proof}
Since ${\mathcal A}=S_{3c}(H/Z\times\Gg/\Gq)$, clearly ${\mathcal A}$
is a dense subspace of ${\mathcal H}$.  The inclusion map (i.\,e.
$\Lambda(f)=f$) carries ${\mathcal A}$ into ${\mathcal H}$.  Now
for $f,g\in{\mathcal A}$,
\begin{align}
&\varphi(g^*\times f)=\int g^*(\tilde{x},\tilde{y},r)f(0-\tilde{x},0-\tilde{y},
r)\bar{e}\bigl[\eta_{\lambda}(r)\beta(\tilde{x},0-\tilde{y})\bigr]\,
d\tilde{x}d\tilde{y}dr  \notag \\
&=\int\overline{g(-\tilde{x},-\tilde{y},r)}\bar{e}\bigl[\eta_{\lambda}(r)
\beta(\tilde{x},\tilde{y})\bigr]f(-\tilde{x},-\tilde{y},r)\bar{e}\bigl
[\eta_{\lambda}(r)\beta(\tilde{x},-\tilde{y})\bigr]\,d\tilde{x}d\tilde{y}dr
\notag \\
&=\int\overline{g(\tilde{x},\tilde{y},r)}f(\tilde{x},\tilde{y},r)
\,d\tilde{x}d\tilde{y}dr=\langle f,g\rangle
=\bigl\langle\Lambda(f),\Lambda(g)\bigr\rangle.  \notag
\end{align}
Consider now the left-multiplication representation $\pi_{\varphi}$.
Then for $f,\xi\in{\mathcal A}$,
\begin{align}
\bigl(\pi_{\varphi}(f)\bigr)\bigl(\Lambda(\xi)\bigr)(x,y,r)&:=
\bigl(\Lambda(f\times\xi)\bigr)(x,y,r)=(f\times\xi)(x,y,r)  \notag \\
&=\int f(\tilde{x},\tilde{y},r)\xi(x-\tilde{x},y-\tilde{y},r)\bar{e}\bigl
[\eta_{\lambda}(r)\beta(\tilde{x},y-\tilde{y})\bigr]\,d\tilde{x}d\tilde{y}
\notag \\
&=L_f\xi(x,y,r).  \notag
\end{align}
This shows that $\pi_{\varphi}$ is just the ${}^*$-representation $L$
of equation \eqref{(representation)}.  Since ${\mathcal A}=\Lambda
({\mathcal A})$ is dense in ${\mathcal H}$, it is also clear that
$\overline{\pi_{\varphi}({\mathcal A}){\mathcal H}}^{\|\ \|_2}
={\mathcal H}$, which means that $\pi_{\varphi}(=L)$ is non-degenerate.
\end{proof}

 We are ready to show that ${\mathcal A}(\subseteq{\mathcal H})$ is a 
``left Hilbert algebra'' (See definition below.).
\begin{defn}(\cite{Cm2}, \cite{St}) 
By a {\em left Hilbert algebra\/}, we mean an involutive algebra
${\mathcal U}$ equipped with a scalar product such that the involution
is an antilinear preclosed mapping in the associated Hilbert space
${\mathcal H}$ and such that the left-multiplication representation
$\pi$ of ${\mathcal U}$ is non-degenerate, bounded, and involutive.
\end{defn}

\begin{prop}\label{lefthilbertA}
The algebra ${\mathcal A}$, together with its inner product inherited
from that of ${\mathcal H}$, is a left Hilbert algebra.
\end{prop}

\begin{proof}
We view ${\mathcal A}=\Lambda({\mathcal A})\subseteq{\mathcal H}$.
It is an involutive algebra equipped with the inner product inherited
from that of ${\mathcal H}$.  Since $\varphi$ is a trace, the map
$f\mapsto f^*$ is not just closable, but it is actually isometry
and hence bounded.  Note that for every $f,g\in{\mathcal A}$,
we have:
$$
\langle f^*,g\rangle=\varphi(g^*\times f^*)=\varphi(f^*\times g^*)=\langle
g^*,f\rangle,
$$
where we used the property that $\varphi$ is a trace.  The remaining
conditions for ${\mathcal A}$ being a left Hilbert algebra are
immediate consequences of the previous proposition.
\end{proof}

\begin{rem}
 The closure of the involution on ${\mathcal A}$ is often denoted by $T$.
The map $T$ is a closed, anti-linear map on ${\mathcal H}$, and ${\mathcal A}$
is a core for $T$ (Actually, $T$ is bounded.).  Define $\nabla=T^*T$.  Clearly,
${\mathcal A}\subseteq D(\nabla)$ and $\nabla(f)=f$ for $f\in{\mathcal A}$.
In other words, $\nabla=\operatorname{Id}$.  The polar decomposition of $T$
is given by $T=J\nabla^{\frac{1}{2}}$, where $\nabla$ is as above and $J$
is an anti-unitary operator.  Obviously in our case, $T=J$.  The ``modular
operator'' $\nabla$ plays an important role in the formulation of the KMS
property.  But as we see here, we can ignore $\nabla$ from now on, all due to
the property that $\varphi$ is a trace.
\end{rem}

 Since we have a left Hilbert algebra structure on ${\mathcal A}$, we
can apply the result of Combes (\cite{Cm2}, \cite{St}) to obtain a weight
extending $\varphi$.  Although it is true that we do not necessarily
have to rely a lot on the theory of weights on $C^*$-algebras (Since
$\varphi$ is a trace in our case, we could use even earlier results of
Dixmier), we nevertheless choose here the more general approach.  The
advantage is that the process will remain essentially the same even
in more difficult examples where we may encounter non-tracial weights.

\begin{theorem}\label{c*weight}
There is a faithful, lower semi-continuous weight on the $C^*$-algebra
$A$ extending the linear functional $\varphi$.  We will use the notation
$\varphi_A$ to denote this weight.
\end{theorem}

\begin{proof}
The representation $\pi_{\varphi}(=L)$ generates the von Neumann algebra
$L({\mathcal A})''$ on the Hilbert space ${\mathcal H}$.  It would be
actually the von Neumann algebra $M_A$ generated by $A$.  On this von
Neumann algebra, we can define as in the below a faithful, semi-finite,
normal weight $\tilde{\varphi}$ (See Theorem 2.11 of \cite{Cm2}.):

For $p\in L({\mathcal A})''$ and $p\ge0$, define $\tilde{\varphi}(p)$ by
$$
\tilde{\varphi}(p)=\left\{\begin{array}{ll}
      \|\xi\|^2=\langle\xi,\xi\rangle  & {\text {if $\exists\xi\in
                {\mathcal A}''$ such that $p^{1/2}=\pi_{\varphi}(\xi)$}} \\
      +\infty  & {\text {otherwise}}
      \end{array}\right .
$$
Here ${\mathcal A}''$ denotes the set of ``left bounded elements''
\cite[\S2]{Cm2}.

We restrict this normal weight to the $C^*$-algebra $\overline
{L({\mathcal A})}^{\|\ \|}$ (norm-closure).  Then the restriction is
a faithful, lower semi-continuous weight.  Since $\pi_{\varphi}(=L)$
extends from ${\mathcal A}$ to an isomorphism $A\cong\overline
{L({\mathcal A})}^{\|\ \|}$, we can use this isomorphism to obtain
the faithful, lower semi-continuous weight (to be denoted by
$\varphi_A$) on $A$.

 It is clear from the construction that $\varphi_A$ extends the linear
functional $\varphi$ on ${\mathcal A}$.  To see this explicitly, suppose
$f\in{\mathcal A}$.  Then $\pi_{\varphi}(f)^*\pi_{\varphi}(f)\in L
({\mathcal A})''$.  According to the theory of left Hilbert algebras,
we then have $\pi_{\varphi}(f)^*\pi_{\varphi}(f)\in{{\mathfrak M}_{\tilde
{\varphi}}}^+$ and
$$
\tilde{\varphi}\bigl(\pi_{\varphi}(f)^*\pi_{\varphi}(f)\bigr)=\bigl\langle
\Lambda(f),\Lambda(f)\bigr\rangle=\langle f,f\rangle=\varphi(f^*f).
$$
But $\pi_{\varphi}(f)^*\pi_{\varphi}(f)=\pi_{\varphi}(f^*f)\in\overline
{L({\mathcal A})}^{\|\ \|}\cong A$, and since $\varphi_A$ is the
restriction of $\tilde{\varphi}$ to $A$, it follows that $\pi_{\varphi}
(f^*f)\in{\mathfrak M}_{\varphi_A}$, and $\varphi_A\bigl(\pi_{\varphi}
(f^*f)\bigr)=\varphi(f^*f)$.  By using polarization, we conclude that
in general, $L({\mathcal A})\subseteq{\mathfrak M}_{\varphi_A}$ and
$$
\varphi_A\bigl(\pi_{\varphi}(f)\bigr)=\varphi(f),\qquad \forall
f\in{\mathcal A}.
$$
\end{proof}

\begin{rem}
>From the proof of the proposition, we can see that $\varphi_A$ is densely
defined (note that we have $L({\mathcal A})\subseteq{\mathfrak M}_{\varphi_A}$).
It is a faithful weight since the linear functional $\varphi$ is faithful
on ${\mathcal A}$.  Since $\varphi_A$ is obtained by restricting the normal
weight $\tilde{\varphi}$ on the von Neumann algebra level, it follows that
it is also KMS (We will not give proof of this here, since $\varphi$ being
a trace makes this last statement redundant: See comment after Definition
\ref{KMS}.).  In the terminology of the first section, $\varphi_A$ is a
``proper'' weight, which is ``faithful'' and ``KMS'' (actually a trace).
\end{rem}

 From now on, let us turn our attention to the weight $\varphi_A$.  Consider
the GNS triple associated with $\varphi_A$, given by the following ingredients:
\begin{itemize}
\item ${\mathcal H}_{\varphi_A}={\mathcal H}$
\item $\Lambda_{\varphi_A}:{\mathfrak N}_{\varphi_A}\to{\mathcal H}$.  The
proof of the previous theorem suggests that for $a\in{\mathfrak N}_{\varphi_A}$,
there exists a unique ``left bounded'' element $v\in{\mathcal H}$.  We define
$\Lambda_{\varphi_A}(a)=v$.
\item $\pi_{\varphi_A}:A\to{\mathcal B}({\mathcal H})$ is the inclusion map.
\end{itemize}
Note that for $f\in{\mathcal A}$, we have: $\Lambda_{\varphi_A}\bigl
(\pi_{\varphi}(f)\bigr)=\Lambda(f)$.  So we know that $\Lambda_{\varphi_A}$
has a dense range in ${\mathcal H}$.

 Define ${\Lambda}_0$ as the closure of the mapping $L({\mathcal A})
\to{\mathcal H}:\pi_{\varphi}(f)\mapsto\Lambda(f)$.  Let us denote by
${\mathcal A}_0$ the domain of ${\Lambda}_0$.  Clearly, $\Lambda_0$ is
a restriction of $\Lambda_{\varphi_A}$.  By using the properties of
$\varphi_A$, including its lower semi-continuity and the ``left invariance''
at the level of the ${}^*$-algebra ${\mathcal A}$, one can improve the left
invariance up to the level of ${\mathcal A}_0$.  One can also show that
${\mathcal A}_0={\mathfrak N}_{\varphi_A}$ and that $L({\mathcal A})$ is
a core for $\Lambda_{\varphi_A}$ (We can more or less follow the discussion
in section 6 of \cite{KuVD}.).  From these results, the left invariance
of $\varphi_A$ can be proved at the $C^*$-algebra level (A similar result
can be found in Corollary 6.14 of \cite{KuVD}.).

 However, we plan to present a somewhat different proof of the left
invariance, which is in the spirit of Van Daele's recently developed method
\cite{VDaxb}, \cite{VDoamp}.  The main strategy is to show that there exists
a faithful, semi-finite, normal weight $\mu$ on ${\mathcal B}({\mathcal H})$
such that at least formally, $\mu(ba)=\varphi_B(b)\varphi_A(a)$ for $b\in B$,
$a\in A$. [See Appendix (Section~6) for the definition of the ``dual'' $B$
and of the weight $\varphi_B$.]

\begin{prop}\label{mu}
Let $\gamma$ be the unbounded operator on ${\mathcal H}$ having $\Lambda
({\mathcal A})$ as a core and is defined by
$$
\gamma\Lambda(f)=\Lambda(\gamma f),\qquad f\in{\mathcal A},
$$
where $\gamma f\in{\mathcal A}$ is such that $\gamma f(x,y,r)=(e^{2\lambda r})^n
f(x,y,r)$.

Now on ${\mathcal B}({\mathcal H})$, we define a linear functional $\mu$ by
$$
\mu:=\operatorname{Tr}(\gamma\,\cdot).
$$
Then $\mu$ is a faithful, semi-finite, normal weight on ${\mathcal B}({\mathcal H})$
such that for $b\in{\mathfrak N}_{\varphi_B}$ and $a\in{\mathfrak N}_{\varphi_A}$,
$$
\mu(b^*a^*ab)=\varphi_B(b^*b)\varphi_A(a^*a).
$$
\end{prop}

\begin{proof}
The operator $\gamma$ is very much related with the ``modular function''
operator, $\tilde{\delta}$, discussed in section 5 (In our case, $\gamma
={\tilde{\delta}}^{-1}$.).  For more precise definition of $\gamma$, see
Definition 2.6 of \cite{VDoamp}.

Let us verify the last statement, at the dense function algebra level of
$b\in\hat{\mathcal A}(\subseteq M_B)$ and $a\in{\mathcal A}(\subseteq M_A)$.
For this, note that by Lemma \ref{haarB} of Appendix and Proposition \ref
{phitrace}, we have:
$$
\varphi_B(b^*b)=\|b\|_2^2,\qquad {\text {and}}\qquad\varphi_A(a^*a)=\|a\|_2^2.
$$
Meanwhile, by equations \eqref{(lambdaf)} and \eqref{(representation)}, we have:
\begin{align}
&(b^*a^*ab)\xi(x,y,r)  \notag \\
&=\int\overline{b(e^{\lambda r}x,e^{\lambda r}y,\tilde{r}-r)}\,
\overline{a(\tilde{x},\tilde{y},\tilde{r})}e\bigl[\eta_{\lambda}(\tilde{r})
\beta(\tilde{x},y)\bigr]a(\hat{x},\hat{y},\tilde{r})\bar{e}\bigl[\eta_{\lambda}
(\tilde{r})\beta(\hat{x},y+\tilde{y}-\hat{y})\bigr]  \notag \\
&\qquad b\bigl(e^{\lambda\hat{r}}(x+\tilde{x}-\hat{x}),e^{\lambda\hat{r}}
(y+\tilde{y}-\hat{y}),\tilde{r}-\hat{r}\bigr)
\xi(x+\tilde{x}-\hat{x},y+\tilde{y}-\hat{y},\hat{r})\,
d\tilde{r}d\tilde{x}d\tilde{y}d\hat{x}d\hat{y}d\hat{r}.  \notag
\end{align}
If we let $(\xi_l)$ be an orthonormal basis in ${\mathcal H}$, we then have:
\begin{align}
\mu(b^*a^*ab)&=\operatorname{Tr}(\gamma b^*a^*ab)=\sum_l\bigl\langle (\gamma
b^*a^*ab)\xi_l,\xi_l\bigr\rangle  \notag \\
&=\sum_l\left(\int (e^{2\lambda r})^n (b^*a^*ab)\xi_l(x,y,r)\overline{\xi_l(x,y,r)}
\,dxdydr\right)  \notag \\
&=\int(e^{2\lambda r})^n\overline{b(e^{\lambda r}x,e^{\lambda r}y,\tilde{r}-r)}\,
\overline{a(\tilde{x},\tilde{y},\tilde{r})}e\bigl[\eta_{\lambda}(\tilde{r})\beta
(\tilde{x},y)\bigr]  \notag \\
&\qquad a(\tilde{x},\tilde{y},\tilde{r})\bar{e}\bigl[\eta_{\lambda}(\tilde{r})
\beta(\tilde{x},y)\bigr]b\bigl(e^{\lambda r}x,e^{\lambda r}y,\tilde{r}-r\bigr)\,
dxdydrd\tilde{x}d\tilde{y}d\tilde{r}  \notag \\
&=\int\overline{b(x,y,r)}\,\overline{a(\tilde{x},\tilde{y},\tilde{r})}
a(\tilde{x},\tilde{y},\tilde{r})b(x,y,r)\,dxdydrd\tilde{x}d\tilde{y}d\tilde{r} 
\notag \\
&=\|b\|_2^2\|a\|_2^2=\varphi_B(b^*b)\varphi_A(a^*a).  \notag
\end{align}
We used the change of variables.

Since $\hat{\mathcal A}$ and ${\mathcal A}$ generate the von Neumann algebras
$M_B$ and $M_A$, while $M_B M_A$ is $\sigma$-strongly dense in ${\mathcal B}
({\mathcal H})$ (see Lemma \ref{MbMa} of Appendix), this will characterize $\mu$.
The properties of $\mu$ being faithful, semi-finite, and normal follow from those
of $\varphi_A$ and $\varphi_B$, as well as the fact that $\mu$ is a trace.
\end{proof}

 The significance of the above proposition is that for a certain well-chosen
element $b\in B$, we know that $A\ni a\mapsto\mu(b^*ab)$ is a scalar multiple
of $\varphi_A$.  This observation is useful in our proof of the left invariance
of $\varphi_A$.  Before we present our main theorem, let us first introduce a lemma
on the linear forms $\omega_{\xi,\eta}$.

\begin{lem}\label{omegalem}
Let $\xi,\eta\in{\mathcal H}$ and consider $\omega_{\xi,\eta}$,
as defined earlier.  If $(\xi_k)$ forms an orthonormal basis of
${\mathcal H}$, we have:
$$
\omega_{\xi,\eta}(ab)=\sum_k\omega_{\xi_k,\eta}(a)\omega_{\xi,\xi_k}(b),
\qquad a,b\in{\mathcal B}({\mathcal H}).
$$
\end{lem}

\begin{proof}
We have:
\begin{align}
\omega_{\xi,\eta}(ab)&=\langle ab\xi,\eta\rangle=\langle b\xi,a^*\eta
\rangle  \notag \\
&=\sum_k\langle b\xi,\xi_k\rangle\langle\xi_k,a^*\eta\rangle
=\sum_k\langle b\xi,\xi_k\rangle\langle a\xi_k,\eta\rangle  \notag \\
&=\sum_k\omega_{\xi,\xi_k}(b)\omega_{\xi_k,\eta}(a).  \notag
\end{align}
\end{proof}

 The following theorem shows the left invariance of $\varphi_A$, as defined
by equation \eqref{(left)}.

\begin{theorem}\label{leftHaar}
For any positive element $a\in A$ such that $\varphi_A(a)<\infty$, and
for $\omega\in A^*_+$, we have:
$$
\varphi_A\bigl((\omega\otimes\operatorname{id})(\Delta a)\bigr)
=\omega(1)\varphi_A(a).
$$
\end{theorem}

\begin{proof}
As stated above, let $a\in{{\mathfrak M}_{\varphi_A}}^{+}$ and let $\omega
\in A^*_+$.  Without loss of generality, we can assume that $\omega$ is
a (positive) vector state.  That is, we can assume that there is a vector
$\zeta\in{\mathcal H}$ such that $\omega=\omega_{\zeta,\zeta}$.

Now consider $(\omega\otimes\operatorname{id})(\Delta a)=(\omega_{\zeta,
\zeta}\otimes\operatorname{id})(\Delta a)$.  For our purposes, it is more
convenient to express $\Delta a$ in terms of the ``dual'' multiplicative
unitary operator defined in Lemma \ref{hatU} in Appendix: From Proposition
\ref{lambdaf}, we know that $\Delta a={\widehat{U_A}}^*(1\otimes a)
\widehat{U_A}$.  If we let $(\xi_k)$ be an orthonormal basis in ${\mathcal H}$,
we would then have:
\begin{align}
(\omega\otimes\operatorname{id})(\Delta a)&=(\omega_{\zeta,\zeta}\otimes
\operatorname{id})\bigl({\widehat{U_A}}^*(1\otimes a)\widehat{U_A}\bigr)   
\notag \\
&=\sum_k\bigl[(\omega_{\xi_k,\zeta}\otimes\operatorname{id})({\widehat{U_A}}^*)
\bigr]a\bigl[(\omega_{\zeta,\xi_k}\otimes\operatorname{id})(\widehat{U_A})\bigr] 
=\sum_k {v_k}^*a^{\frac{1}{2}}a^{\frac{1}{2}}v_k.   \notag
\end{align}
The sum is convergent in the $\sigma$-weak topology on the von Neumann
algebra $M_A$ (Use Lemma \ref{omegalem}.).  For convenience, we let
$v_k=(\omega_{\zeta,\xi_k}\otimes\operatorname{id})(\widehat{U_A})\bigl(\in
{\mathcal B}({\mathcal H})\bigr)$.  Note that since $\widehat{U_A}$ is
unitary, the operators $v_k$ have the property that for the orthonormal
basis $(\xi_l)$ of ${\mathcal H}$, we have:
\begin{align}\label{(vk)}
\sum_k\langle v_k\xi_l,v_k\xi_j\rangle&=\bigl\langle\widehat{U_A}(\zeta
\otimes\xi_l),\widehat{U_A}(\zeta\otimes\xi_j)\bigr\rangle   \notag \\
&=\langle\zeta\otimes\xi_l,\zeta\otimes\xi_j\rangle=\langle\zeta,\zeta
\rangle\langle\xi_l,\xi_j\rangle.
\end{align}

Next, suggested by Proposition \ref{mu} and the comments following it,
let us choose a fixed element $b\in\hat{\mathcal A}(\subseteq
{\mathfrak N}_{\varphi_B})$, so that we have:
$$
\varphi_A(a)=\left(\frac{1}{\|b\|_2^2}\right)\mu(b^*ab),\quad {\text { for
$a\in{\mathfrak M}_{\varphi_A}$}}.
$$
Then combining these observations, we have the following:
\begin{align}
&\varphi_A\bigl((\omega_{\zeta,\zeta}\otimes\operatorname{id})(\Delta a)
\bigr)=\varphi_A\left(\sum_k {v_k}^*a^{\frac{1}{2}}a^{\frac{1}{2}}v_k\right)
=\sum_k\varphi_A({v_k}^*a^{\frac{1}{2}}a^{\frac{1}{2}}v_k)    \notag \\
&=\frac{1}{\|b\|_2^2}\sum_k\mu(b^*{v_k}^*a^{\frac{1}{2}}a^{\frac{1}{2}}v_kb)
=\frac{1}{\|b\|_2^2}\sum_{k,l}\operatorname{Tr}(\gamma b^*{v_k}^*
a^{\frac{1}{2}}a^{\frac{1}{2}}v_kb\xi_l,\xi_l)  \notag \\
&=\frac{1}{\|b\|_2^2}\sum_{k,l}\langle v_k\gamma^{\frac{1}{2}}
a^{\frac{1}{2}}b\xi_l,v_k\gamma^{\frac{1}{2}}a^{\frac{1}{2}}b\xi_l\rangle
\notag \\
&=\frac{1}{\|b\|_2^2}\sum_l\langle\zeta,\zeta\rangle\langle
\gamma^{\frac{1}{2}}a^{\frac{1}{2}}b\xi_l,\gamma^{\frac{1}{2}}
a^{\frac{1}{2}}b\xi_l\rangle\qquad\quad {\text { by equation \eqref{(vk)}}}
\notag \\
&=\frac{1}{\|b\|_2^2}\langle\zeta,\zeta\rangle\mu(b^*ab)=\langle\zeta,
\zeta\rangle\varphi_A(a)=\|\omega\|\varphi_A(a)=\omega(1)\varphi_A(a).
\notag
\end{align}
\end{proof}

 As we remarked in section 1, proving this ``weak'' version of the left
invariance is enough.  In this way, we have shown that $\varphi_A$ is a proper,
faithful, KMS (tracial) weight on $(A,\Delta)$, which is left invariant.
This satisfies the requirement of Definition \ref{lcqgdefn}.

 We now need to talk about the right invariant weight on $(A,\Delta)$.
Again by viewing $(A,\Delta)$ as a ``quantized $C_0(G)$'', we try
to build the weight from the right Haar measure of $G$ (The group
structure of $G$ as defined in Definition 1.6 of \cite{BJKp2}
immediately gives us the natural choice for its right Haar measure.).
Just as we did at the beginning of this section, this suggestion
lets us to consider the linear functional $\psi$ on ${\mathcal A}$,
as described below.

\begin{defn}
On ${\mathcal A}$, define a linear functional $\psi$ by
$$
\psi(f)=\int f(0,0,r)(e^{-2\lambda r})^n\,dr.
$$
\end{defn}

 It is helpful to realize that at the level of the ${}^*$-algebra
${\mathcal A}$, we have: $\psi=\varphi\circ S$, where $S$ is the
antipodal map we defined in Proposition \ref{antipode}.  Indeed,
for $f\in{\mathcal A}$, we have:
\begin{align}
\varphi\bigl(S(f)\bigr)&=\int\bigl(S(f)\bigr)(0,0,r)\,dr=\int(e^{2\lambda r})^n
\bar{e}\bigl[\eta_{\lambda}(r)\beta(0,0)\bigr]f(0,0,-r)\,dr  \notag \\
&=\int f(0,0,r)(e^{-2\lambda r})^n\,dr=\psi(f).  \notag
\end{align}
Therefore, to extend $\psi$ to the $C^*$-algebra level, we may consider
$\psi_A=\varphi_A\circ S$, where $S$ is now regarded as an antiautomorphism
on $A$.

\begin{rem}
Defining $\psi_A=\varphi_A\circ S$ is not entirely correct: In general,
the ``antipode'' $S$ may not be defined everywhere and can be unbounded.
However, even in the general case, the antipode can be always written
in the form $S=R\tau_{-\frac{i}{2}}$ (``polar decomposition'' of $S$),
where $\tau$ is the so-called ``scaling group'' and $R$ is the ``unitary
antipode''.  In our case, $\tau\equiv\operatorname{Id}$ and $R=S$ (See
section 4.).  The correct way of defining $\psi_A$ would be: $\psi_A=
\varphi_A\circ R$, which is true in general.
\end{rem}

 Since $R$ is an (anti-)automorphism on $A$, it follows that $\psi_A
=\varphi_A\circ R$  is clearly a faithful, lower semi-continuous, densely
defined KMS weight on $A$, extending $\psi$.  Checking the ``right
invariance'' is straightforward, if we use the property of $R$.

\begin{theorem}\label{rightHaar}
Let $\psi_A=\varphi_A\circ R$.  It is a proper, faithful, KMS (and tracial)
weight on $A$.  It is also ``right invariant''.  That is, for $a\in
{{\mathfrak M}_{\psi_A}}^+$ and for $\omega\in A^*_+$, we have:
$$
\psi_A\bigl((\operatorname{id}\otimes\omega)(\Delta a)\bigr)
=\omega(1)\psi_A(a).
$$
\end{theorem}

\begin{proof}
Recall from Proposition \ref{antipode} that $R(=S)$ satisfies
$(R\otimes R)(\Delta a)=\chi\bigl(\Delta(R(a))\bigr)$, where $\chi$
denotes the flip.  We thus have:
\begin{align}
\psi_A\bigl((\operatorname{id}\otimes\omega)(\Delta a)\bigr)
&=\varphi_A\bigl(R((\operatorname{id}\otimes\omega)(\Delta a))\bigr)
\notag \\
&=\varphi_A\bigl((\operatorname{id}\otimes\omega)((R\otimes R)(\Delta a))
\bigr)  \notag \\
&=\varphi_A\bigl((\operatorname{id}\otimes\omega)(\chi(\Delta(R(a))))\bigr)
=\varphi_A\bigl((\omega\otimes\operatorname{id})(\Delta(R(a)))\bigr) 
\notag \\
&=\omega(1)\varphi_A\bigl(R(a)\bigr)\qquad\qquad\qquad\qquad{\text
{$\varphi_A$: left invariant}}  \notag \\
&=\omega(1)\psi_A(a).  \notag
\end{align}
\end{proof}

 We thus have the weight $\psi_A$ on $(A,\Delta)$, satisfying the
requirement of Definition \ref{lcqgdefn}.  For another characterization
of the right invariant weight, see section 5.

 Finally, we are now able to say that $(A,\Delta)$ is indeed a
{\em ($C^*$-algebraic) locally compact quantum group\/}, in the
sense of \cite{KuVa}.

\begin{theorem}\label{main}
The pair $(A,\Delta)$, together with the weights $\varphi_A$ and
$\psi_A$ on it, is a $C^*$-algebraic locally compact quantum group.
\end{theorem}

\begin{proof}
Combine the results of Proposition \ref{comultiplication} and
Proposition \ref{density} on the comultiplication $\Delta$.
Theorem \ref{c*weight} and Theorem \ref{leftHaar} gives the
left invariant weight $\phi_A$, while Theorem \ref{rightHaar}
gives the right invariant weight $\psi_A$.  By Definition \ref
{lcqgdefn}, we conclude that $(A,\Delta)$ is a {\em (reduced)
$C^*$-algebraic quantum group\/}.
\end{proof}

\section{Antipode}

 According to the general theory (by Kustermans and Vaes \cite{KuVa}),
the result of Theorem \ref{main} is enough to establish our main
goal of showing that $(A,\Delta)$ is indeed a $C^*$-algebraic locally
compact quantum group (satisfying Definition \ref{lcqgdefn}).

 Assuming both the left invariant and the right invariant weights in
the definition may look somewhat peculiar, while there is no mention
on the antipode.  However, using these rather simple set of axioms,
Kustermans and Vaes could prove additional properties for $(A,\Delta)$,
so that it can be legitimately called a locally compact quantum group.
They first construct a manageable multiplicative unitary operator (in
the sense of \cite{BS} and \cite{Wr7}) associated with $(A,\Delta)$.
[In our case, this unitary operator $W$ coincides with our $\widehat{U_A}$
defined in Appendix.]  More significantly, they then construct the
antipode and its polar decomposition.  The uniqueness (up to scalar
multiplication) of the Haar weight is also obtained.

 An aspect of note through all this is that in this new definition,
the ``left (or right) invariance'' of a weight has been formulated
without invoking the antipode, while a characterization of the antipode
is given without explicitly referring to any invariant weights.
This is much simpler and is a fundamental improvement over earlier
frameworks, where one usually requires certain conditions of the type:
$$
(\operatorname{id}\otimes\varphi)\bigl((1\otimes a)(\Delta b)\bigr)
=S\bigl((\operatorname{id}\otimes\varphi)((\Delta a)(1\otimes b))\bigr).
$$
It is also more natural.  Note that in the cases of ordinary locally
compact groups or Hopf algebras in the purely algebraic setting, the
 axioms of the antipode do not have to require any relationships to
invariant measures.

 For details on the general theory, we will refer the reader to \cite
{KuVa}.  What we plan to do in this section is to match the general
theory with our specific example.  Let us see if we can re-construct $S$
from $(A,\Delta)$.

 By general theory (\cite{Wr7}, \cite{KuVa}), the antipode, $S$, can
be characterized such that $\bigl\{(\omega\otimes\operatorname{id})
(U_A):\omega\in{\mathcal B}({\mathcal H})_*\bigr\}$ is a core for $S$
and
\begin{equation}\label{(antipodeS)}
S\bigl((\omega\otimes\operatorname{id})(U_A)\bigr)=(\omega\otimes
\operatorname{id})({U_A}^*),\qquad\omega\in{\mathcal B}({\mathcal H})_*.
\end{equation}
It is a closed linear operator on $A$.  The domain $D(S)$ is a
subalgebra of $A$ and $S$ is antimultiplicative: i.\,e. $S(ab)
=S(b)S(a)$, for any $a,b\in D(S)$.  The image $S\bigl(D(S)\bigr)$
coincides with $D(S)^*$ and $S\bigl(S(a)^*\bigr)^*=a$ for any
$a\in D(S)$.  The operator $S$ admits the (unique) ``polar
decomposition'': $S=R\tau_{-\frac{i}{2}}$, where $R$ is the
``unitary antipode'' and $\tau_{-\frac{i}{2}}$ is the analytic
generator of a certain one parameter group $(\tau_t)_{t\in\mathbb{R}}$
of ${}^*$-automorphisms of $A$ (called the ``scaling group'').

\begin{rem}
In \cite{KuVa}, the scaling group and the unitary antipode are
constructed first (using only the multiplicative unitary operator
and the invariant weights), from which they define the antipode
via $S=R\tau_{-\frac{i}{2}}$.  The characterization given above
is due to Woronowicz \cite{Wr7}.
\end{rem}

 To compare $S$ given by equation \eqref{(antipodeS)} with our
own $S$ defined in Proposition \ref{antipode}, let us again consider
$\omega_{\xi,\eta}$.  From the proof of Proposition \ref{A(U)}, we
know that
$$
(\omega_{\xi,\eta}\otimes\operatorname{id})(U_A)=L_f,
$$
where
$$
f(\tilde{x},\tilde{y},r)=\int\xi(\tilde{x},\tilde{y},\tilde{r}+r)
(e^{\lambda r})^n\overline{\eta(e^{\lambda r}\tilde{x},e^{\lambda r}
\tilde{y},\tilde{r})}\,d\tilde{r}.
$$

 We can carry out a similar computation for $(\omega_{\xi,\eta}
\otimes\operatorname{id})({U_A}^*)$.  For $\zeta\in{\mathcal H}$,
we would have (again using change of variables):
\begin{align}
\bigl(S(L_f)\bigr)\zeta(x,y,r)&=S\bigl((\omega_{\xi,\eta}\otimes
\operatorname{id})(U_A)\bigr)\zeta(x,y,r) \notag \\
&=(\omega_{\xi,\eta}\otimes\operatorname{id})({U_A}^*)\zeta(x,y,r) 
\notag \\
&=\int\bigl({U_A}^*(\xi\otimes\zeta)\bigr)(\tilde{x},\tilde{y},
\tilde{r};x,y,r)\overline{\eta(\tilde{x},\tilde{y},\tilde{r})}
\,d\tilde{x}d\tilde{y}d\tilde{r}  \notag \\
&=\int g(\tilde{x},\tilde{y},r)\zeta(x-\tilde{x},y-\tilde{y},r)
\bar{e}\bigl[\eta_{\lambda}(r)\beta(\tilde{x},y-\tilde{y})
\bigr]\,d\tilde{x}d\tilde{y}  \notag  \\
&=L_g\zeta(x,y,r),  \notag
\end{align}
where
$$
g(\tilde{x},\tilde{y},r)=\int(e^{\lambda r})^n\xi(-e^{\lambda r}
\tilde{x},-e^{\lambda r}\tilde{y},\tilde{r}-r)\overline{\eta(-\tilde
{x},-\tilde{y},\tilde{r})}\bar{e}\bigl[\eta_{\lambda}(r)\beta(\tilde{x},
\tilde{y})\bigr]\,d\tilde{r}.
$$
This means that $S(f)=g$.  Comparing with $f$, we see that
$$
\bigl(S(f)\bigr)(x,y,r)=g(x,y,r)=(e^{2\lambda r})^n f(-e^{\lambda r}
\tilde{x},-e^{\lambda r}\tilde{y},-r)\bar{e}\bigl[\eta_{\lambda}(r)
\beta(x,y)\bigr].
$$
This is exactly the expression we gave in Proposition \ref{antipode},
verifying that our situation agrees perfectly with the general theory.

 Since we have already seen that $S:A\to A$ is an antiautomorphism
(defined everywhere on $A$), the uniqueness of the polar decomposition
implies that $R=S$ and $\tau\equiv\operatorname{Id}$.

 As a final comment on the general theory, we point out that after
one defines the antipode as $S=R\tau_{-\frac{1}{2}}$, one proves
that for $a,b\in{\mathfrak N}_{\psi_A}$,
$$
S\bigl((\psi_A\otimes\operatorname{id})((a^*\otimes1)(\Delta b))\bigr)
=(\psi_A\otimes\operatorname{id})\bigl(\Delta(a^*)(b\otimes1)\bigr).
$$
In this way, one can ``define'' $S$, as well as give a stronger
version of the invariance of $\psi_A$.  The fact that this result
could be obtained from the defining axioms (as opposed to being
one of the axioms itself) was the significant achievement of
\cite{KuVa}.

\section{Modular function}

 To motivate the modular function of $(A,\Delta)$, let us re-visit
our right invariant weight $\psi_A$.  We will keep the notation
of Section 3.  Recall that at the level of the dense ${}^*$-algebra
${\mathcal A}$, the right invariant weight is given by the linear
functional $\psi$:
$$
\psi(f)=\int f(0,0,r)(e^{-2\lambda r})^n\,dr.
$$

 Let us consider the Hilbert space ${\mathcal H}_R$, which will
be the GNS Hilbert space for $\psi$.  It is defined such that
${\mathcal H}_R={\mathcal H}$ as a space and the inner product
on it is defined by
$$
\langle f,g\rangle_R=\int f(x,y,r)\overline{g(x,y,r)}
(e^{-2\lambda r})^n\,dr.
$$
Let $\Lambda_R$ be the inclusion map $\Lambda_R:{\mathcal A}
\hookrightarrow{\mathcal H}_R$.  We can see easily that for
$f,g\in{\mathcal A}$,
$$
\bigl\langle\Lambda_R(f),\Lambda_R(g)\bigr\rangle_R=\bigl\langle
\Lambda(f),\Lambda(\delta g)\bigr\rangle.
$$
Here $\delta g\in{\mathcal A}$ defined by $\delta g(x,y,r)=
(e^{-2\lambda r})^n g(x,y,r)$.

 For motivational purposes, let us be less rigorous for the
time being.  Observe that working purely formally, we can
regard $\delta g$ as follows:
\begin{align}
\delta g(x,y,r)&=(e^{-2\lambda r})^n g(x,y,r)  \notag \\
&=\int\delta(\tilde{x},\tilde{y},r)g(x-\tilde{x},y-\tilde{y},r)
\bar{e}\bigl[\eta_{\lambda}(r)\beta(\tilde{x},y-\tilde{y})\bigr]
\,d\tilde{x}d\tilde{y}  \notag \\
&=(\delta\times g)(x,y,r),  \notag
\end{align}
where $\delta$ is considered as a (Dirac delta type) ``function''
in the $(x,y,r)$ variables such that
\begin{align}
&\delta(x,y,r)=0,\qquad\qquad\qquad{\text {(if $x\ne0$ or $y\ne0$)}}
\notag \\
&\int\delta(x,y,r)\,dxdy=(e^{-2\lambda r})^n.   \notag
\end{align}
At the level of the functions in the $(x,y,z)$ variables, $\delta$
corresponds to the following ``function'' (we may use partial
Fourier transform, again purely formally):
$$
\delta(x,y,z)=\int e\bigl[(-e^{-\lambda r}p)\cdot x+(-e^{-\lambda r}q)
\cdot y+(-r)z\bigr]\,dpdqdr.
$$
In this formulation, we see an indication of the inverse operation
on $G$ (Note that in $G$, we have $(p,q,r)^{-1}=(-e^{-\lambda r}p,
-e^{-\lambda r}q,-r)$.).

 These remarks modestly justifies our intention to call $\delta$
a ``modular function''.  Certainly, we see that $\delta$ plays an
important role relating $\langle\ ,\ \rangle_R$ and $\langle\ ,\ 
\rangle$, or in other words, relating $\psi$ and $\varphi$.  A word of
caution is that $\delta$ is not bounded and not exactly a function.
What we plan to do here is to make this notion precise in the
$C^*$-algebra setting.

 Note that since ${\mathcal A}$ is a dense subspace of ${\mathcal H}$,
we may already regard the map $\tilde{\delta}:{\mathcal A}\ni g\mapsto
\delta g\in{\mathcal A}$ as an operator on ${\mathcal H}$.  It would be
an (unbounded) operator affiliated with the von Neumann algebra $L(A)''
=M_A$, since for an arbitrary element $b\in L(A)'$ and $g\in{\mathcal A}
(\subseteq A)$, we would have:
$$
b\tilde{\delta}g=\tilde{\delta}gb=\tilde{\delta}bg.
$$
By viewing ${\mathcal A}$ as a dense subspace of ${\mathcal H}$, we
conclude that $\tilde{\delta}$ commutes with $b\in L(A)'$, proving
our claim that $\tilde{\delta}$ is affiliated with $M_A$.

 We may pull down the operator $\tilde{\delta}$ to the $C^*$-algebra
level, and obtain an operator affiliated with $A$, in the $C^*$-algebra
setting (c.\,f. in the sense of Woronowicz \cite{Wr4}).  So define first
a closed linear (unbounded) operator, $N$, from ${\mathcal H}$ into
${\mathcal H}_R$ such that $\Lambda({\mathcal A})$ is a core for $N$
and
$$
N\Lambda(f):=\Lambda_R(f),\qquad f\in{\mathcal A}.
$$
Then $N$ is a densely defined, injective operator with dense range.
Note also that $\bigl\langle N\Lambda(f),\Lambda_R(g)\bigr\rangle_R
=\bigl\langle\Lambda(f),\Lambda(\delta g)\bigr\rangle$.  So we have
$\Lambda_R({\mathcal A})\subseteq D(N^*)$, and
$$
N^*\Lambda_R(g)=\Lambda(\delta g),\qquad g\in{\mathcal A}.
$$
Consider the following operator (which will be the ``modular function'').
Clearly, $\Lambda({\mathcal A})\subseteq D(\delta_A)$ and $\delta_A
\Lambda(f)=\Lambda(\delta f)$.

\begin{defn}
Define $\delta_A=N^*N$.  It is an injective, positive (unbounded)
operator on ${\mathcal H}$.
\end{defn}

 By general theory, we can say that $\delta_A$ is the appropriate
definition of the ``modular function'' in the $C^*$-algebra setting.

\begin{theorem}
Let $\delta_A$ be defined as above.  Then the following properties
hold.
\begin{enumerate}
\item $\delta_A$ is an operator affiliated with the $C^*$-algebra $A$.
\item $\Delta(\delta_A)=\delta_A\otimes\delta_A$.
\item $\tau_t(\delta_A)=\delta_A$ and $R(\delta_A)=\delta_A^{-1}$.
\item $\psi(a)=\varphi(\delta_A^{\frac{1}{2}}a\delta_A^{\frac{1}{2}})$,
for $a\in{\mathcal A}$.
\end{enumerate}
\end{theorem}

\begin{proof}
It is not difficult to see that $\delta_A$ is cut down from
the operator $\tilde{\delta}$.  For proof of the statements,
see \cite[\S7]{KuVa} or see \cite[\S8]{KuVD}.  There are also
important relations relating the modular automorphism groups
corresponding to $\varphi_A$ and $\psi_A$, but in our case they
become trivial.
\end{proof}

\section{Appendix: An alternative formulation of the dual}

 The aim of this Appendix is to present a dual counterpart to our locally
compact quantum group $(A,\Delta)$, which is slightly different (though
isomorphic) from $(\hat{A},\hat{\Delta})$ defined in Section 2.  It would be
actually the Hopf $C^*$-algebra having the opposite multiplication and
the opposite comultiplication to $(\hat{A},\hat{\Delta})$.  To avoid a
lengthy discussion, we plan to give only a brief treatment.  But we will
include results that are relevant to our main theorem in Section 3.

 Let us define $(B,\Delta_B)$, by again using the language of multiplicative
unitary operators.  We begin with a lemma, which is motivated by the general
theory of multiplicative unitary operators \cite{BS}.

\begin{lem}\label{hatU}
Let $j\in{\mathcal B}({\mathcal H})$ be defined by
$$
j\xi(x,y,r)=(e^{\lambda r})^n\bar{e}\bigl[\eta_{\lambda}(r)\beta(x,y)\bigr]
\xi(-e^{\lambda r}x,-e^{\lambda r}y,-r).
$$
Then $j$ is a unitary operator such that $j^2=1$.  Moreover, the operator
$\widehat{U_A}$ defined by
$$
\widehat{U_A}=\Sigma(j\otimes1)U_A(j\otimes1)\Sigma,\qquad
{\text {$\Sigma$ denotes the flip}}
$$
is multiplicative unitary and is regular.  For $\xi\in{\mathcal H}$, we
specifically have:
\begin{align}
\widehat{U_A}\xi(x,y,r,x',y',r')&=e\bigl[\eta_{\lambda}(r)\beta(e^{\lambda(r'-r)}x',
y)\bigr] \notag \\
&\quad\xi(x+e^{\lambda(r'-r)}x',y+e^{\lambda(r'-r)}y',r;x',y',r'-r). \notag
\end{align}
\end{lem}

\begin{rem}
The proof is straightforward.  What is really going on is that the triple
$({\mathcal H},U_A,j)$ forms a {\em Kac system\/}, in the terminology of Baaj
and Skandalis (See section 6 of \cite{BS}.).  The operator $j$ may be written
as $j=\hat{J}J=J\hat{J}$, where $\hat{J}$ is the anti-unitary operator defined
in the proof of Proposition \ref{antipode}, while $J$ is the anti-unitary
operator determining the ${}^*$-operation of $A$ as mentioned in the remark
following Proposition \ref{lefthilbertA}.  We will have an occasion to say
more about these operators in our future paper.
\end{rem}

\begin{defn}\label{dualopp}
Let $\widehat{U_A}$ be the multiplicative unitary operator obtained above.  Define
$(B,\Delta_B)$ as follows:
\begin{enumerate}
\item Let ${\mathcal A}(\widehat{U_A})$  be the subspace of ${\mathcal B}({\mathcal H})$
defined by
$$
{\mathcal A}(\widehat{U_A})=\bigl\{(\omega\otimes\operatorname{id})(\hat{U_A}):
\omega\in{\mathcal B}({\mathcal H})_*\bigr\}.
$$
Then ${\mathcal A}(\widehat{U_A})$ is a subalgebra of ${\mathcal B}({\mathcal H})$,
and the subspace ${\mathcal A}(\widehat{U_A}){\mathcal H}$ forms a total set in
${\mathcal H}$.
\item The norm-closure in ${\mathcal B}({\mathcal H})$ of the algebra ${\mathcal A}
(\widehat{U_A})$ is the $C^*$-algebra $B$.  The $\sigma$-strong closure of ${\mathcal A}
(\widehat{U_A})$ in ${\mathcal B}({\mathcal H})$ will be the von Neumann algebra $M_B$.
\item For $b\in{\mathcal A}(\widehat{U_A})$, define ${\Delta_B}(b)$ by ${\Delta_B}(b)=
\widehat{U_A}(b\otimes 1){\widehat{U_A}}^*$.  Then ${\Delta_B}$ can be extended to the
comultiplication on $B$, and also to the level of the von Neumann algebra $M_B$.
\end{enumerate}
\end{defn}

 We do not give the proof here, since it is essentially the same as in Propositions
\ref{A(U)} and \ref{dualqg}.  Let us just add a brief clarification: By a
comultiplication on $B$, we mean a non-degenerate $C^*$-homomorphism ${\Delta_B}:
B\to M(B\otimes B)$ satisfying the coassociativity; whereas by a comultiplication
on $M_B$, we mean a unital normal ${}^*$-homomorphism ${\Delta_B}:M_B\to M_B\otimes
M_B$ satisfying the coassociativity.  In the following proposition, we give a more
specific description of the $C^*$-algebra $B$.

\begin{prop}\label{lambdaf}
\begin{enumerate}
\item For $\omega\in{\mathcal B}({\mathcal H})_*$, we let ${\lambda}(\omega)
=(\omega\otimes\operatorname{id})(\widehat{U_A})$.  Then we have: ${\lambda}(\omega)
=j{\rho}(\omega)j$, where ${\rho}(\omega)=(\operatorname{id}\otimes\omega)(U_A)\in
\hat{\mathcal A}(U_A)$ as defined in equation \eqref{(rhof)}.
\item Let $\hat{\mathcal A}$ be the space of Schwartz functions in the $(x,y,r)$
variables having compact support in the $r$ variable.  For $f\in\hat{\mathcal A}$,
define the operator ${\lambda}_f\in{\mathcal B}({\mathcal H})$ by
\begin{equation}\label{(lambdaf)}
({\lambda}_f\zeta)(x,y,r)=\int f(e^{\lambda\tilde{r}}x,e^{\lambda\tilde{r}}y,r-\tilde{r})
\zeta(x,y,\tilde{r})\,d\tilde{r}. 
\end{equation}
Then the $C^*$-algebra $B$ is generated by the operators ${\lambda}_f$.  By partial
Fourier transform, we can also show that $B\cong\lambda\bigl(C^*(G)\bigr)$, where
$\lambda$ is the left regular representation of $C^*(G)$.
\item For any $f,g\in\hat{\mathcal A}$, we have: $[{\rho}_f,{\lambda}_g]=0$.
Actually, $M_B={M_{\hat{A}}}'$, where $M_{\hat{A}}$ is the von Neumann algebra
generated by $\hat{A}$.
\item For $b\in(\hat{A},\hat{\Delta})$, we have: $({\lambda}\otimes{\lambda})
(\hat{\Delta}b)=\widehat{U_A}\bigl({\lambda}(b)\otimes 1\bigr){\widehat{U_A}}^*$.
\item Dually, there exists an alternative characterization of the $C^*$-algebra
$A$:
$$
A=\overline{\bigl\{(\operatorname{id}\otimes\omega)(\widehat{U_A}):\omega\in
{\mathcal B}({\mathcal H})_*\bigr\}}^{\|\ \|}.
$$
And for $a\in A$, we have: $(L\otimes L)(\Delta a)={\widehat{U_A}}^*\bigl(1\otimes
L(a)\bigr)\widehat{U_A}$, where $L$ is the regular representation of $A$ defined
in Section 2.
\end{enumerate}
\end{prop}

\begin{proof}
See Proposition 6.8 of \cite{BS}.  For instance, for the first statement, note that:
$$
{\lambda}(\omega)=(\operatorname{id}\otimes\omega)\bigl((j\otimes1)U_A(j\otimes1)\bigr)
=j{\rho}(\omega)j.
$$
The second statement is a consequence of this result.  We can also give a direct proof,
just as in Propositions \ref{A(U)} and \ref{dualqg}.  Actually, we have:
$$
j{\rho}_fj={\lambda}_{\tilde{f}},\qquad f\in\hat{\mathcal A},
$$
where $\tilde{f}(x,y,r)=\bar{e}\bigl[\eta_{\lambda}(r)\beta(x,y)\bigr]
f(-e^{\lambda r}x,-e^{\lambda r}y,-r)$.

For the third and fourth statements, we can again refer to general theory (Proposition
6.8 of \cite{BS}), or we can just give a direct proof.  Since we see (up to partial
Fourier transform in the $x$ and $y$ variables) that ${\rho}_f$ and ${\lambda}_f$
are essentially the right and left regular representations of $C^*(G)$, the result
follows easily.  The last statement is also straightforward (similar to Proposition
\ref{A(U)}).
\end{proof}

\begin{rem}
The above proposition implies that at least at the level of the dense
subalgebra of functions, $B$ has an opposite algebra structure to that
of $\hat{A}$.  Meanwhile, (4) above implies that $(\hat{A},\hat{\Delta})
\cong(B,\Delta_B)$ as Hopf $C^*$-algebras.
\end{rem}

 It turns out that working with $(B,{\Delta}_B)$ and $M_B={M_{\hat{A}}}'$ is more
convenient in our proof of Theorem \ref{leftHaar}.  Here are a couple of lemmas that
are useful in Section 3.  Similar results exist for $\hat{A}$ and $M_{\hat{A}}$.
We took light versions of the proofs.

\begin{lem}\label{MbMa}
Let $M_A$ and $M_B$ be the enveloping von Neumann algebras of $A$ and $B$.  We have:
\begin{enumerate}
\item $\widehat{U_A}\in M_A\otimes M_B$.
\item $M_B\cap M_A=\mathbb{C}1$.
\item The linear space $M_B M_A$ is $\sigma$-strongly dense in ${\mathcal B}
({\mathcal H})$.
\end{enumerate}
\end{lem}

\begin{proof}
The first statement is immediate from general theory, once we realize
(see previous proposition) that $\widehat{U_A}$ determines $B$ and $A$
(as well as $M_B$ and $M_A$).  We also have: $\widehat{U_A}\in M(A\otimes B)$.
The remaining two results also follow from the same realization.  For
instance, we could modify the proof of Proposition 2.5 of \cite{VDoamp}.
\end{proof}

\begin{lem}\label{haarB}
On $\hat{\mathcal A}\subseteq B$, consider a linear functional $\varphi_B$
defined by
$$
\varphi_B({\lambda}_f)=\int f(x,y,0)\,dxdy.
$$
It can be extended to a faithful, semi-finite, normal weight ${\tilde
{\varphi}}_B$ on $M_B$.
\end{lem}

\begin{rem}
The idea for proof of this lemma is pretty much the same as the early part of section 3
(but $\varphi_B$ is no longer a trace).  It turns out that $\varphi_B$ will be an
(invariant) Haar weight for $(B,\Delta_B)$, although for our current purposes, this
result is not immediately necessary.  We will make all these clear in our future paper.
Meanwhile, by a straightforward calculation using equation \eqref{(lambdaf)}, we see
easily that:
$$
\varphi_B({{\lambda}_f}^*{\lambda}_f)=\|f\|_2^2.
$$
This last result will be useful in the proof of Proposition \ref{mu}.
\end{rem}

 Using multiplicative unitary operators, we can also give notions that are
analogous to ``opposite dual'' or ``co-opposite dual'' \cite{BJKqhg}.  For
a more careful discussion on the duality, as well as on the notion of quantum
double, refer to our future paper.



\bibliographystyle{amsplain}

\providecommand{\bysame}{\leavevmode\hbox to3em{\hrulefill}\thinspace}

\end{document}